%% file: main.tex
\documentclass{arxivpaper}
\usepackage{manuscript}
\renewcommand{\fullversionref}[2]{#1}

\makeatletter
\jmlrproceedings{}{}
\jmlrvolume{}
\jmlryear{}
\jmlrpages{}
\jmlrsubmitted{}
\jmlrpublished{}
\jmlrworkshop{}
\renewcommand*{\@titlefoot}{\hfill\small\rmfamily\thepage\hfill}
\makeatother

\begin{document}
\maketitle
\input{body}
\input{endmatter}
\appendix
\input{appendix}
\end{document}

%% file: body.tex
\begin{abstract}
  We propose CBO-MPPI, a sampling-based predictive controller that applies consensus-based optimization (CBO) to the free energy underlying model predictive path integral control (MPPI). Under hard constraints, this objective rewards sampling neighborhoods containing feasible, low-cost trajectories under perturbation. Standard MPPI performs a local stochastic search by updating one nominal control sequence from Gaussian perturbations around it. In CBO-MPPI, we estimate the free energy at multiple nominal centers from local rollouts. These estimates define a shared consensus that guides drift and diffusion, and the final selected center supplies an MPPI weighted-mean output. For the associated mean-field model, we establish well-posedness and, under additional assumptions, convergence toward a free-energy minimizer. In simulations of obstacle navigation, pendulum swing-up, and planar-arm manipulation under equal rollout budgets, at high temperatures, CBO-MPPI attains lower free energy than MPPI and multi-start MPPI and maintains greater obstacle clearance than CBO-MPC.
\end{abstract}


\section{Introduction}
\label{sec:introduction}

Model predictive control (MPC) repeatedly solves a finite-horizon optimization problem from the measured state. Its ability to incorporate nonlinear dynamics and constraints directly has led to widespread use in robotics, autonomous driving, and process control \citep{Grune2017Nonlinear}. The resulting optimization problems can be high-dimensional and nonconvex, with nonsmooth or discontinuous objectives arising from obstacle avoidance and hard constraints. Solving these problems efficiently within the available computation time is a central challenge in MPC.

Model predictive path integral control (MPPI) offers a sampling-based approach to the MPC subproblem \citep{Kappen2005Path,Theodorou2010Generalized,Williams2017Model,Williams2018InformationTheoretic}. It generates perturbed control sequences around a nominal sequence and updates that sequence by exponentially weighting the rollout costs. The update requires rollout costs without differentiating the dynamics or cost, and the rollouts can be evaluated in parallel. The objective underlying this update is a free energy defined through an exponential average of perturbed trajectory costs. Under hard constraints, a low free energy identifies a sampling neighborhood containing substantial probability mass on feasible, low-cost trajectories. With fixed Gaussian covariance, one MPPI update approximates a preconditioned gradient step on this free-energy objective \citep{Fazlyab2026Model}. Because MPPI repeatedly updates one nominal sequence, it has a local search structure and may miss distant basins with lower free energy.

We propose \emph{CBO-MPPI} to address this local-search limitation by evolving multiple MPPI sampling centers with consensus-based optimization (CBO). CBO is a gradient-free particle method for global optimization, driven by consensus and random exploration \citep{Pinnau2017consensusbased,Carrillo2021consensusbased,Fornasier2024ConsensusBased}. In CBO-MPPI, the particles are the centers, and the consensus favors those with lower free-energy estimates from local rollouts. All centers move toward this common consensus instead of taking independent MPPI steps. Since drift and diffusion depend on the distance to the consensus rather than on the local landscape, a center at a local minimum keeps moving as long as it differs from the consensus. At the final round, the center with the lowest estimated free energy supplies the control through the usual MPPI weighted mean.

Several related approaches also aim to improve exploration in sampling-based optimal control. Closest to our work, CBO has been applied directly to trajectory costs in finite-horizon optimal control problems \citep{BorghiMonSep0100:00:00UTC2025Model,Sun2026Consensusbased}. Other methods improve sampling efficiency by adapting where control sequences are sampled through adaptive importance sampling \citep{Asmar2023Model} and learned proposal distributions \citep{Sacks2023Learning}. Variational and Stein MPC represent alternative solutions through multimodal control distributions \citep{Okada2020Variational,Lambert2021Stein}, while multi-modal MPPI uses maneuver-specific proposals \citep{Bertipaglia2025MultiModal}. SVG-MPPI guides MPPI toward a target mode by updating guide particles with surrogate gradients \citep{Honda2024Steina}. CBO-MPPI differs in how its centers interact: a free-energy-weighted consensus governs their gradient-free drift and diffusion.

In this paper, we make three contributions. First, CBO-MPPI applies CBO to the MPPI free energy over a population of sampling centers. Second, Gaussian smoothing yields an exact, unprojected mean-field model under compactness and bounded-measurability assumptions, even when hard constraints make the trajectory cost discontinuous. Additional assumptions yield finite-accuracy global convergence. Third, we compare CBO-MPPI with MPPI, multi-start MPPI, and CBO-MPC under equal rollout budgets on Navigation, Many Obstacles, Pendulum, and Planar Arm benchmarks. At high temperatures, CBO-MPPI has the lowest free energy among MPPI controllers and greater obstacle clearance than CBO-MPC; its MPPI free-energy advantage persists across tested settings.

\section{Problem Formulation and Method}
\label{sec:problem_setting}

We define the constrained trajectory cost and MPPI free energy, then present MPPI and CBO-MPPI.

\subsection{Constrained Trajectory Cost and Free Energy}

We consider a discrete-time system with state and input dimensions $n$ and $m$, respectively, and state-transition map $f:\R^n\times\R^m\to\R^n$.
At each control period, we fix the current state $x\in\R^n$ and a prediction horizon of $H\geq1$ steps. A candidate control sequence $U=(u_0,\ldots,u_{H-1})\in\R^D$, with $u_t\in\R^m$ and $D=mH$, generates the predicted trajectory from $x_0=x$ through
\begin{equation}
  x_{t+1}=f(x_t,u_t),\quad t=0,\ldots,H-1.
  \label{eq:discrete_time_dynamics}
\end{equation}
For a stage cost $\ell:\R^n\times\R^m\to\R$ and terminal cost $g:\R^n\to\R$, the trajectory cost along~\eqref{eq:discrete_time_dynamics} is
\begin{equation}
  J_x(U)=\sum_{t=0}^{H-1}\ell(x_t,u_t)+g(x_H).
  \label{eq:open_loop_trajectory_cost}
\end{equation}
For the control problem, we specify a state-constraint set $\mathcal C_x\subseteq\R^n$ and an actuator box $\mathcal C_u\subset\R^m$. Let $\mathcal C_u^H:=\prod_{t=0}^{H-1}\mathcal C_u$ denote the $H$-fold Cartesian product of the actuator set. The admissible input sequences satisfy both the input bounds and the state constraints along the predicted trajectory in~\eqref{eq:discrete_time_dynamics}:
\[
  \mathcal C=\{U\in\mathcal C_u^H:x_t\in\mathcal C_x,\quad t=1,\ldots,H\}.
\]
We define the extended cost by $\bar J_x(U)=J_x(U)$ for $U\in\mathcal C$ and $\bar J_x(U)=+\infty$ otherwise. At each control period, MPC approximately minimizes the cost in~\eqref{eq:open_loop_trajectory_cost} over $U\in\mathcal C$ and applies the first input of the resulting sequence to the plant. We denote this first input by $[U]_0=u_0$.

For this fixed control problem, we make the following assumptions.

\begin{assumption}[Admissible control set]
  \label{ass:admissible_set}
  The admissible control set $\mathcal C\subset\mathbb{R}^D$ is compact and has positive Lebesgue measure.
\end{assumption}

\begin{assumption}[Open-loop cost]
  \label{ass:cost}
  For the fixed current state $x\in\R^n$, the open-loop cost $J_x$ in~\eqref{eq:open_loop_trajectory_cost}, restricted to $\mathcal C$, is Borel measurable and bounded: there exist constants $J_{\min},J_{\max}\in\mathbb{R}$ with $J_{\min}\leq J_x(U)\leq J_{\max}$ for all $U\in\mathcal C$.
\end{assumption}

To define the free energy, fix a temperature $\tau>0$ and a sampling covariance $\Sigma\in\mathbb{R}^{D\times D}$ with $\Sigma\succ0$. We use $\mathcal N(0,\Sigma)$ for the Gaussian probability distribution and $\varphi_\Sigma$ for its Lebesgue density.
For a nominal sequence $M\in\R^D$ and $\epsilon\sim\mathcal N(0,\Sigma)$, we define the free energy by
\begin{equation}
  F_x(M):=-\tau\log\E\!\left[e^{-\bar J_x(M+\epsilon)/\tau}\right]
  =-\tau\log\!\int_{\mathcal C}\varphi_\Sigma(U-M)e^{-J_x(U)/\tau}\,dU.
  \label{eq:mppi_free_energy}
\end{equation}
Under Assumptions~\ref{ass:admissible_set} and \ref{ass:cost}, the expectation is strictly positive and finite, and Gaussian convolution gives $F_x\in C^\infty(\R^D)$. \fullversionref{Appendix~\ref{subsec:proof_free_energy_properties}}{The supplementary material} proves these properties.

The free energy can be decomposed by conditioning on the perturbed sequence being admissible:
\begin{equation}
  F_x(M)=-\tau\log\Pr(M+\epsilon\in\mathcal C)
  -\tau\log\E\!\left[e^{-J_x(M+\epsilon)/\tau}\mid M+\epsilon\in\mathcal C\right].
  \label{eq:hard_constraint_free_energy_decomposition}
\end{equation}
The first term penalizes a small probability of drawing an admissible sequence; the second is a soft minimum of the costs conditional on admissibility. Thus a low value of $F_x(M)$ identifies a sampling neighborhood around $M$ that contains sufficiently many admissible perturbations with low trajectory cost. The temperature $\tau$ controls how strongly the exponential average favors lower costs and scales the feasibility penalty. This interpretation motivates optimizing free energy under hard constraints.

\subsection{MPPI as Single-Center Search}

MPPI draws $L$ samples $U_j=M+\epsilon_j$ around a nominal sequence $M$, where $\epsilon_j\sim\mathcal N(0,\Sigma)$, and evaluates the rollout costs $S_j=\bar J_x(U_j)$ along their predicted trajectories \citep{Williams2017Model,Williams2018InformationTheoretic,Fazlyab2026Model}. We use the batch average $\widehat Z=L^{-1}\sum_{j=1}^L e^{-S_j/\tau}$ to normalize weights $w_j$, $j=1,\ldots,L$, that favor lower-cost samples. We update the nominal sequence to the weighted mean of the samples:
\begin{equation}
  w_j=\frac{e^{-S_j/\tau}}{L\widehat Z},
  \qquad M^+=\sum_{j=1}^{L}w_jU_j.
  \label{eq:sample_mppi_update}
\end{equation}
The first input $[M^+]_0$ supplies the current control action, and the remaining sequence provides an initial guess for the next control period.

We define the exact-expectation MPPI update as the exponentially weighted mean
\begin{equation}
  \Phi(M):=
  \frac{\int_{\R^D}U\varphi_\Sigma(U-M)e^{-\bar J_x(U)/\tau}\,dU}
       {\int_{\R^D}\varphi_\Sigma(U-M)e^{-\bar J_x(U)/\tau}\,dU}
  =M-\tau^{-1}\Sigma\nabla F_x(M).
  \label{eq:mppi_gradient_identity}
\end{equation}
The second equality is the fixed-covariance Gaussian MPPI gradient identity \citep{Fazlyab2026Model}. Under Assumptions~\ref{ass:admissible_set} and \ref{ass:cost}, boundedness of $J_x$ on compact $\mathcal C$ provides an integrable bound for the weighted Gaussian-kernel derivatives, justifying differentiation under the integral. The finite-batch update~\eqref{eq:sample_mppi_update} approximates this mean. Thus the expectation update $\Phi(M)$ is a gradient step on $F_x$ with preconditioner $\Sigma/\tau$. Each MPPI update thus estimates a descent direction from samples around a single center, so the search remains local. \fullversionref{Appendix~\ref{subsec:proof_free_energy_properties} states these identities and verifies their applicability, and Algorithm~\ref{alg:mppi} gives the control-period flow used by MPPI}{The supplementary material states these identities, verifies their applicability, and gives the control-period flow used by MPPI}.

\subsection{CBO-MPPI}
\label{subsec:paper_cbo_mppi}

CBO-MPPI combines local MPPI rollout evaluation with CBO population search on free energy \citep{Pinnau2017consensusbased,Carrillo2021consensusbased}. We maintain $N$ nominal centers $M_i^k\in\R^D$, $i=1,\ldots,N$, as particles and perform $K$ evaluation rounds, indexed by $k=0,\ldots,K-1$. As in single-center MPPI, $M$ denotes a sampling center and $U$ a candidate input sequence. At round $k$, each center uses the MPPI sampling construction: we draw $L$ sequences $U_{i,j}^k=M_i^k+\epsilon_{i,j}^k$, with $\epsilon_{i,j}^k\sim\mathcal N(0,\Sigma)$, and evaluate $S_{i,j}^k=\bar J_x(U_{i,j}^k)$ from the current state $x$. The batch average used to normalize the MPPI weights in~\eqref{eq:sample_mppi_update} also yields an estimate of the free energy in~\eqref{eq:mppi_free_energy} at each center:
\begin{equation}
  \widehat Z_i^k=\frac1L\sum_{j=1}^{L}e^{-S_{i,j}^k/\tau},
  \qquad \widehat F_i^k=-\tau\log\widehat Z_i^k.
  \label{eq:finite_sample_free_energy}
\end{equation}
For a CBO inverse temperature $\alpha>0$ and a round with at least one finite free-energy estimate, we define the center weights and consensus point
\begin{equation}
  \omega_i^k=\frac{\exp(-\alpha\widehat F_i^k)}
  {\sum_{r=1}^N\exp(-\alpha\widehat F_r^k)},
  \qquad \widehat M_\alpha^k:=\sum_{i=1}^N\omega_i^k M_i^k.
  \label{eq:finite_sample_consensus}
\end{equation}
After each of the first $K-1$ evaluation rounds, the shared point $\widehat M_\alpha^k$ drives the consensus drift and diffusion of all centers. The Euler--Maruyama scheme gives
\begin{equation}
  \begin{aligned}
    M_i^{k+1}
    ={}&M_i^k
    -\lambda\Delta s\left(M_i^k-\widehat M_\alpha^k\right)
    +\sigma_{\mathrm{cbo}}
    \left\|M_i^k-\widehat M_\alpha^k\right\|
    \xi_i^k\sqrt{\Delta s},
    \qquad
    \xi_i^k\sim\mathcal N(0,I_D).
  \end{aligned}
  \label{eq:discrete_cbo_mppi_update}
\end{equation}
Here $\lambda>0$ is the drift strength, $\sigma_{\rm cbo}>0$ the diffusion strength, and $\Delta s>0$ the optimization step size; $\|\cdot\|$ denotes the Euclidean norm and $I_D$ the $D$-dimensional identity matrix. The independent Gaussian vectors $\xi_i^k$ produce diffusion proportional to each center's distance from consensus.

At the final round, we select $i^\star\in\argmin_i\widehat F_i^{K-1}$ and reuse its batch for the MPPI weights and output
\begin{equation}
  w_j^\star=
  \frac{e^{-S_{i^\star,j}^{K-1}/\tau}}
       {L\widehat Z_{i^\star}^{K-1}},
  \qquad U_{\rm out}=\sum_{j=1}^{L}w_j^\star U_{i^\star,j}^{K-1},
  \quad u_{\rm apply}=[U_{\rm out}]_0.
  \label{eq:weighted_mean_output}
\end{equation}
We warm-start the next MPC period by shifting the final evaluated centers and adding independent smooth jitter. The operator $\operatorname{Shift}$ removes the first input, shifts the remaining $H-1$ inputs one step to the left, and appends a terminal input; the numerical implementation appends zero. Viewing each center as an $H\times m$ array, we use a basis $B_{\rm jit}\in\R^{H\times R_{\rm jit}}$ of $R_{\rm jit}$ low-frequency cosine modes along the horizon and set $\zeta_i=\sigma_{\rm jit}B_{\rm jit}A_i$. Here $\sigma_{\rm jit}$ sets the jitter amplitude, and each $A_i\in\R^{R_{\rm jit}\times m}$ has independent standard Gaussian entries.
Algorithm~\ref{alg:cbo_mppi_free_energy} gives the control-period procedure, and \fullversionref{Appendix~\ref{sec:paper_experimental_details}}{the supplementary material} specifies the shift and jitter construction shared by the population methods.

\begin{algorithm}[h]
  \algocfconts{alg:cbo_mppi_free_energy}{%
    \caption{CBO-MPPI: population search on MPPI free energy}%
  }{%
    \begin{algorithmic}[1]
      \REQUIRE Current state $x$, horizon $H$, initialized or warm-started centers $M_i^0\in\R^{mH}$, $i=1,\ldots,N$; local rollouts per center $L$, evaluation rounds $K$; $\tau,\alpha>0$, $\Sigma\succ0$; $\lambda,\sigma_{\rm cbo},\Delta s>0$; $\sigma_{\rm jit}$ and $B_{\rm jit}$
      \ENSURE Applied input $u_{\rm apply}$ and next-period centers $M_i^{0,\rm next}$
      \FOR{$k=0,\ldots,K-1$}
      \FOR{all centers $i=1,\ldots,N$ in parallel}
      \STATE Draw $\epsilon_{i,j}^k\sim\mathcal N(0,\Sigma)$ and set $U_{i,j}^k=M_i^k+\epsilon_{i,j}^k$, $j=1,\ldots,L$
      \STATE Roll out from $x$ and evaluate $S_{i,j}^k=\bar J_x(U_{i,j}^k)$ for all $j$ in parallel
      \STATE Compute $\widehat F_i^k$ by~\eqref{eq:finite_sample_free_energy}
      \ENDFOR
      \STATE Compute the consensus $\widehat M_\alpha^k$ by~\eqref{eq:finite_sample_consensus}
      \IF{$k<K-1$}
      \STATE Draw $\xi_i^k\sim\mathcal N(0,I_D)$ and update every center by~\eqref{eq:discrete_cbo_mppi_update}
      \ENDIF
      \ENDFOR
      \STATE Select $i^\star\in\argmin_i\widehat F_i^{K-1}$
      \STATE Reuse the selected batch to compute $U_{\rm out}$ by~\eqref{eq:weighted_mean_output}; apply $u_{\rm apply}=[U_{\rm out}]_0$
      \FOR{all $i=1,\ldots,N$}
      \STATE Draw $A_i$ with independent $\mathcal N(0,1)$ entries and set $\zeta_i=\sigma_{\rm jit}B_{\rm jit}A_i$
      \STATE Set $M_i^{0,\rm next}=\operatorname{Shift}(M_i^{K-1})+\zeta_i$
      \ENDFOR
    \end{algorithmic}
  }
\end{algorithm}

\begin{remark}[Temperature parameters]
For a round with at least one finite free-energy estimate, suppressing the index $k$, the center weights in~\eqref{eq:finite_sample_consensus} satisfy $\omega_i\propto e^{-\alpha\widehat F_i}=\widehat Z_i^{\alpha\tau}$. Thus $\tau$ controls exponential weighting of rollouts within each batch, while $\alpha\tau$ is the exponent applied to the batch averages. Since $\widehat Z_i$ itself depends on $\tau$, fixing $\alpha\tau$ does not make the weights independent of $\tau$.
\end{remark}

\begin{remark}[Input projection and zero-weight batches]
For numerical rollout evaluation, each input sequence is projected onto $\mathcal C_u^H$ before simulating its trajectory. Denoting this projection by $\Pi_{\mathcal C_u^H}$, Algorithm~\ref{alg:cbo_mppi_free_energy} evaluates $\bar J_x(\Pi_{\mathcal C_u^H}(U_{i,j}))$. Rollouts that still violate the state constraints, including obstacles, receive cost $+\infty$ and weight zero. The implemented objective is therefore $F_x^\Pi(M) := -\tau\log\mathbb E\big[e^{-\bar J_x(\Pi_{\mathcal C_u^H}(M+\epsilon))/\tau}\big]$, which~\eqref{eq:finite_sample_free_energy} estimates at each center. All weights use the projected rollout costs, while weighted means and consensus points are formed from the unprojected samples and centers. The applied first input is projected onto $\mathcal C_u$ before being applied to the plant. For zero-weight cases, we use two fallback rules. If every rollout in a local batch is infeasible, we set its free-energy estimate to $+\infty$ and its local output to the projected center. If all center weights vanish, we use the mean of the projected centers as the consensus. \fullversionref{Appendix~\ref{sec:paper_experimental_details}}{The supplementary material} gives the common numerical protocol.
\end{remark}

\begingroup
\submissionlayoutfalse
\section{Theoretical Analysis}
\label{sec:theoretical_analysis}

We analyze the mean-field CBO-MPPI dynamics for a fixed open-loop problem.
We establish well-posedness using the regularity of the free energy induced by Gaussian smoothing \citep{Fazlyab2026Model}.
The regularity and growth bounds verified here provide the basis for applying existing CBO analyses. In particular, \fullversionref{Appendix~\ref{subsec:paper_global_convergence}}{the supplementary material} develops a finite-accuracy global-convergence result under additional minimizer, initialization, and drift--diffusion assumptions.

\label{subsec:ideal_mean_field_model}
We define the mean-field CBO-MPPI model for a fixed current state $x\in\R^n$ and $H,\tau>0,\Sigma\succ0$, with the associated admissible set $\mathcal C$. The model uses the exact, unprojected free energy $F_x$ in~\eqref{eq:mppi_free_energy}.
Let $\mathcal P(\mathbb R^D)$ denote the set of Borel probability measures on $\mathbb R^D$; for $1\leq p<\infty$, let $\mathcal P_p(\mathbb R^D)$ be its finite-$p$th-moment subset, endowed with the Wasserstein-$p$ distance $W_p$.
The driving process $B=(B_s)_{s\geq0}$ is a standard $D$-dimensional Brownian motion.
The mean-field model is the McKean--Vlasov SDE
\begin{equation}
  \ifsubmissionlayout
    \begin{aligned}
      dM_s={}   & -\lambda(M_s-m_\alpha[\rho_s])\,ds \\
                & +\sigma_{\mathrm{cbo}}
      \|M_s-m_\alpha[\rho_s]\|\,dB_s,                \\
      \rho_s={} & \operatorname{Law}(M_s),
    \end{aligned}
  \else
    dM_s=-\lambda(M_s-m_\alpha[\rho_s])\,ds
    +\sigma_{\mathrm{cbo}}
    \|M_s-m_\alpha[\rho_s]\|\,dB_s,
    \qquad
    \rho_s=\operatorname{Law}(M_s),
  \fi
  \label{eq:mean_field_cbo_mppi}
\end{equation}
The law $\rho_s$ represents the distribution of nominal centers at optimization time $s$.
Here $m_\alpha[\rho]$ is the free-energy-weighted consensus point for a particle distribution $\rho\in\mathcal P(\mathbb R^D)$, defined by
\begin{equation}
  m_\alpha[\rho]
  :=
  \frac{
    \int_{\mathbb R^D}M e^{-\alpha F_x(M)}\rho(dM)
  }{
    \int_{\mathbb R^D}e^{-\alpha F_x(M)}\rho(dM)
  }.
  \label{eq:mean_field_consensus}
\end{equation}
Under Assumptions~\ref{ass:admissible_set} and \ref{ass:cost}, $F_x$ is continuous and bounded below, so $m_\alpha[\rho]\in\mathbb R^D$ is well defined for every $\rho\in\mathcal P_1(\mathbb R^D)$.
The following theorem establishes well-posedness, using regularity and growth bounds for $F_x$ proved in \fullversionref{Appendix~\ref{sec:proof}}{the supplementary material}.

\begin{samepage}
\begin{theorem}[Well-posedness]
  \label{prop:mean_field_well_posedness}
  Suppose Assumptions~\ref{ass:admissible_set} and \ref{ass:cost} hold, let $\rho_0\in\mathcal P_4(\mathbb R^D)$, and let $\lambda,\sigma_{\mathrm{cbo}},\alpha>0$. For every $S>0$, there exists a unique nonlinear strong process $M$ with continuous sample paths satisfying \eqref{eq:mean_field_cbo_mppi} and $\operatorname{Law}(M_0)=\rho_0$. It satisfies $\sup_{0\leq s\leq S}\mathbb E\|M_s\|^4<\infty$, and $\rho_s:=\operatorname{Law}(M_s)$ belongs to $C([0,S],\mathcal P_4(\mathbb R^D))$.
\end{theorem}
\end{samepage}

For the proof, Gaussian-kernel differentiation yields derivatives of $F_x$ without requiring $J_x$ itself to be differentiable. Compactness of $\mathcal C$ gives global regularity and weighted-Lipschitz bounds, while Gaussian tails give quadratic growth. After translating a global minimizer to the origin, these bounds allow us to apply \cite[Theorem~3.2]{Carrillo2016analytical} for existence, uniqueness, and the fourth-moment bound. Jensen's inequality and the Burkholder--Davis--Gundy inequality then establish continuity of the law in $W_4$. The complete proof is in \fullversionref{Appendix~\ref{subsec:proof_mean_field_well_posedness}}{the supplementary material}, where we also establish finite-accuracy global convergence\fullversionref{ (Theorem~\ref{thm:mean_field_global_convergence})}{} under additional assumptions.

\endgroup
\section{Numerical Examples}
\label{sec:numerical_example}
We study four closed-loop tasks: Navigation, Many Obstacles, Pendulum, and Planar Arm. Navigation is the main route-selection example; Many Obstacles uses the same point mass in a cluttered geometry. Pendulum examines swing-up under actuator disturbances, and Planar Arm considers redundant manipulation with hard link-obstacle constraints.

\subsection{Experimental Setting and Metrics}

We compare CBO-MPPI with three baselines. MPPI performs one weighted-mean update per control period. Multi-start MPPI (MS-MPPI) independently applies a local MPPI update to each center after each of the first $K-1$ evaluation rounds. As in CBO-MPPI, its final output uses the MPPI weighted mean at the center with the smallest estimated free energy, and the final evaluated centers are shifted and jittered for the next control period. CBO-MPC treats each candidate control sequence as a particle, evaluates it with one rollout, applies CBO directly to trajectory cost \citep{BorghiMonSep0100:00:00UTC2025Model,Sun2026Consensusbased}, and uses the first input of the lowest-cost particle. \fullversionref{Appendix~\ref{full:sec:proposed_method}}{The supplementary material} gives pseudocode for the three baseline controllers.

We set the controller parameters as follows. We compare the methods at a fixed per-step rollout budget $L_{\mathrm{tot}}$. MPPI evaluates $L$ rollouts in one update, so $L_{\mathrm{tot}}=L$. CBO-MPC evaluates $N$ particles in each of $K$ rounds, with one rollout per particle, so $L_{\mathrm{tot}}=NK$. MS-MPPI and CBO-MPPI evaluate $L$ local rollouts around each of $N$ centers in each of $K$ rounds, so $L_{\mathrm{tot}}=KNL$. Using preliminary timing tests, we selected budgets to fit each task's control period: $L_{\mathrm{tot}}=131{,}072$ for Navigation, Many Obstacles, and Planar Arm and $L_{\mathrm{tot}}=16{,}384$ for Pendulum. MS-MPPI matches CBO-MPPI in $(N,L,K)$, initialization, local perturbations, and warm-start jitter. For each comparison, the MPPI-based controllers share the temperature $\tau$ and sampling covariance $\Sigma$. \fullversionref{Appendices~\ref{sec:paper_experimental_details}--\ref{subsec:paper_arm_details} detail}{The supplementary material details} the controller parameter settings.

For closed-loop evaluation, we define the realized cost $J_{\mathrm{CL}}$ and the free energy $F_{\mathrm{CL}}$ of the realized input sequence. All methods use the same $T$-step evaluation window within each task. For realized states $x_0,\ldots,x_T$ and applied inputs $U_{\mathrm{CL}}:=(u_0,\ldots,u_{T-1})\in\mathbb R^{mT}$, the realized cost is $J_{\mathrm{CL}}:=\sum_{k=0}^{T-1}\ell(x_{k+1},u_k)+g(x_T)$. For $L_{\mathrm e}$ perturbed sequences $U_{\mathrm{CL}}+\epsilon_j$ with zero-mean Gaussian $\epsilon_j\in\mathbb R^{mT}$, the corresponding free energy is
\begin{equation}
  F_{\mathrm{CL}}
  :=-\tau\log\Biggl[
    \frac{1}{L_{\mathrm e}}\sum_{j=1}^{L_{\mathrm e}}
    \exp\!\left(
    -\frac{\bar J_{x_0}\!\left(\Pi_{\mathcal C_u^T}(U_{\mathrm{CL}}+\epsilon_j)\right)}{\tau}
    \right)\Biggr].
  \label{eq:paper_closed_loop_metrics}
\end{equation}
We evaluate $F_{\rm CL}$ after the control episode, with $\bar J_{x_0}$ defined over this window. \fullversionref{Appendix~\ref{sec:paper_experimental_details}}{The supplementary material} specifies the task-dependent perturbation settings and evaluation sample counts.

All controllers and evaluation routines are implemented in Python using JAX/JAXlib~0.8.3 and NumPy~2.4.6. The numerical experiments used Apple M4 Max and NVIDIA RTX 6000 Ada hardware; implementation timings were measured on the former, a MacBook Pro with 128~GB unified memory. Each $t_{\rm ctrl}$ is the median of 15 run-level control-cycle 95th percentiles, with the initial JIT compilation excluded. Code and reproducibility materials are available at \url{https://github.com/inody/cbo_mppi}.

\subsection{Navigation: Objective Geometry and Closed-Loop Search}
\label{subsec:two_obstacles}

This task steers a planar point mass with position $p_k\in\R^2$, velocity $v_k\in\R^2$, and acceleration input $u_k\in\R^2$ through a narrow passage between two obstacles to a target. The dynamics are $p_{k+1}=p_k+\Delta t\,v_k$ and $v_{k+1}=v_k+\Delta t\,u_k$, with time step $\Delta t=0.1$. The initial state is $(p_0,v_0)=(-4,0,0,0)$, and the target position is $(4,0)$. Two obstacles of radius $0.5$ are centered at $(0,\pm0.85)$. We encode obstacle avoidance in $\mathcal C_x$ by excluding states whose positions lie inside either obstacle. The cost penalizes goal-position error, velocity, control effort, and input increments, with terminal penalties on position error and velocity. The narrow central passage and the longer detours are competing routes.

\begin{figure}[t]
  \centering
  \includegraphics[width=\linewidth]{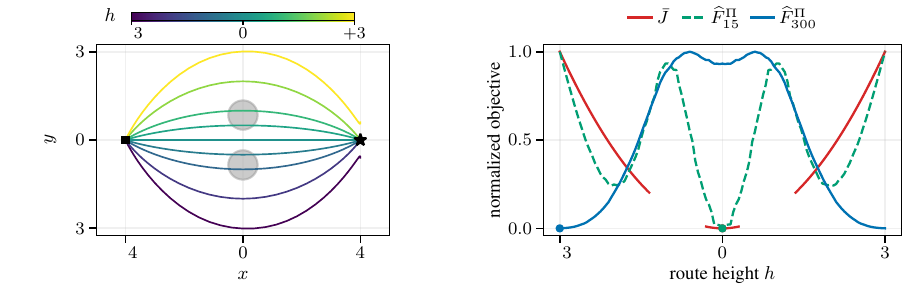}
  \caption{Navigation objective geometry. Left: prescribed open-loop routes colored by signed route height $h$. Right: nominal cost and estimated projected free energy along the same family, with each curve mapped to $[0,1]$ using its own minimum and maximum. In the legend, $\widehat F^\Pi_{15}$ and $\widehat F^\Pi_{300}$ denote the free energy at $\tau=15$ and $\tau=300$.}
  \label{fig:paper_landscapes}
\end{figure}

\paragraph{Route family and objective preference.}
We first compare objectives on artificially constructed open-loop routes before examining trajectories produced by feedback control. The routes are parameterized by signed crossing height $h$: $h=0$ passes through the central gap; positive and negative $h$ bend the routes upward and downward. Figure~\ref{fig:paper_landscapes} (left) shows representative routes. For each route's input sequence, we evaluate the nominal trajectory cost and the projected free energy at $\tau=15$ and $300$. Figure~\ref{fig:paper_landscapes} (right) plots these values against $h$. Perturbations around the narrow center route are more likely to encounter obstacles than those around wide detours, increasing its free energy. The nominal cost and low-temperature free energy favor the center route. The high-temperature free-energy minimum shifts to a detour. \fullversionref{Appendix~\ref{subsec:paper_navigation_details}}{The supplementary material} gives the task specifications, route parameterization, and unnormalized evaluation.

\begin{figure}[t]
  \centering
  \includegraphics[width=\linewidth]{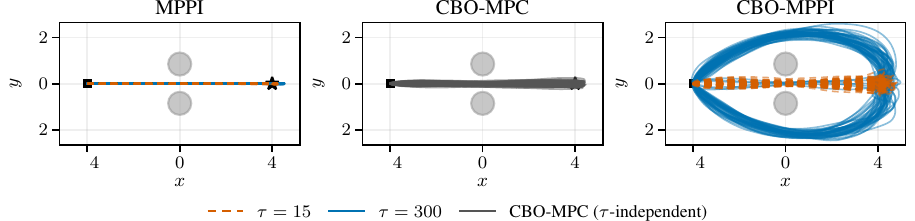}
  \caption{Navigation trajectories for MPPI, CBO-MPC, and CBO-MPPI at equal rollout budget, with all 100 trials per condition.}
  \label{fig:paper_navigation}
\end{figure}

\paragraph{Closed-loop search.}
Table~\ref{tab:paper_main_results} compares median $J_{\rm CL}$ and $F_{\rm CL}$ for all four methods at both temperatures. Figure~\ref{fig:paper_navigation} shows trajectories for MPPI, CBO-MPC, and CBO-MPPI. At $\tau=15$, the three plotted methods primarily use the center passage. At $\tau=300$, MPPI uses that passage in all 100 trials, whereas CBO-MPPI detours in 94/100 trials and has greater median minimum obstacle clearance ($0.717$ versus $0.353$). CBO-MPPI's lower $F_{\rm CL}$ at this temperature accompanies realized cost close to MPPI's. MS-MPPI also lowers $F_{\rm CL}$ relative to MPPI at high temperature and incurs a larger $J_{\rm CL}$; at low temperature, it improves neither cost. CBO-MPPI has lower $F_{\rm CL}$ than MS-MPPI at both temperatures. CBO-MPC has the lowest $J_{\rm CL}$ at both temperatures and the lowest $F_{\rm CL}$ at $\tau=15$; CBO-MPPI has the lowest $F_{\rm CL}$ at $\tau=300$.

CBO-MPPI's route selection follows the temperature-dependent preference of the open-loop landscape; MPPI remains in the center passage. CBO-MPC also obtains smaller $F_{\rm CL}$ with center trajectories by accumulating less position-error cost and achieving higher perturbation feasibility than MPPI. The open-loop landscape therefore helps explain CBO-MPPI's route selection, and $F_{\rm CL}$ additionally reflects closed-loop trajectory cost and perturbation feasibility.

\subsection{Complementary Tasks: Many Obstacles, Pendulum, and Planar Arm}

We conduct experiments on three additional tasks and summarize the four-task comparison in Table~\ref{tab:paper_main_results}. At the selected budgets, CBO-MPPI's $t_{\rm ctrl}$ is below the control period on all four tasks and is comparable to baseline values. Across the obstacle-avoidance tasks, the high-temperature results show larger median minimum clearance for CBO-MPPI than for CBO-MPC: $0.717$ versus $0.319$ on Navigation, $0.458$ versus $0.080$ on Many Obstacles, and $0.132$ versus $0.003$ on Planar Arm. \fullversionref{Appendices~\ref{subsec:paper_navigation_details}--\ref{subsec:paper_arm_details} give}{The supplementary material gives} the detailed setups and complete results for all four tasks.

\paragraph{Many Obstacles.}
This task steers the same point mass through a cluttered obstacle field. At $\tau=15$, CBO-MPC has the lowest median $F_{\rm CL}$, followed by MPPI, CBO-MPPI, and MS-MPPI. At $\tau=3000$, CBO-MPPI has the lowest $F_{\rm CL}$ and MS-MPPI is second; both remain below MPPI and CBO-MPC. CBO-MPC and MPPI attain smaller $J_{\rm CL}$ than CBO-MPPI at high temperature, while MS-MPPI has the largest realized cost. Thus CBO-MPPI's high-temperature free-energy advantage persists in the cluttered setting, with higher realized cost than MPPI and CBO-MPC.

\paragraph{Pendulum.}
This task swings a pendulum from downward to upright under bounded inputs and shared actuator disturbances. CBO-MPPI has the lowest median $F_{\rm CL}$ at both $\tau=1$ and $10$. MS-MPPI has the highest $F_{\rm CL}$ at low temperature. At high temperature, it improves on MPPI and remains above CBO-MPC and CBO-MPPI. At $\tau=1$, CBO-MPPI's smaller free energy comes with larger $J_{\rm CL}$ than MPPI and CBO-MPC. At $\tau=10$, CBO-MPPI lowers both costs relative to MPPI and MS-MPPI; CBO-MPC still has smaller $J_{\rm CL}$. These Pendulum trials extend the observed free-energy improvement beyond obstacle-avoidance tasks.

\paragraph{Planar Arm.}
This task moves a three-link arm around obstacles to a target end-effector position. At $\tau=3$, MS-MPPI has the lowest median $F_{\rm CL}$ and $J_{\rm CL}$; CBO-MPC is second, followed by CBO-MPPI and MPPI. At $\tau=300$, CBO-MPPI attains smaller $F_{\rm CL}$ than MS-MPPI, and both remain below CBO-MPC and MPPI. CBO-MPC has the smallest $J_{\rm CL}$ at this temperature.

\begin{table}[t]
  \centering
  \caption{Median $J_{\rm CL}$ and $F_{\rm CL}$ over 100 trials per task at equal rollout budgets. MS-MPPI uses persistent independent centers; CBO-MPC is fixed across temperatures, with only $F_{\rm CL}$ reevaluated at each $\tau$.}
  \label{tab:paper_main_results}
  \small
  \setlength{\tabcolsep}{3.5pt}
  \begin{tabular}{@{}lr*{4}{rr}@{}}
    \toprule
                   &        & \multicolumn{2}{c}{MPPI} & \multicolumn{2}{c}{MS-MPPI} & \multicolumn{2}{c}{CBO-MPC} & \multicolumn{2}{c}{CBO-MPPI} \\
    \cmidrule(lr){3-4}\cmidrule(lr){5-6}\cmidrule(lr){7-8}\cmidrule(l){9-10}
    Task           & $\tau$ & $J_{\rm CL}$ & $F_{\rm CL}$ & $J_{\rm CL}$ & $F_{\rm CL}$ & $J_{\rm CL}$ & $F_{\rm CL}$ & $J_{\rm CL}$ & $F_{\rm CL}$ \\
    \midrule
    Navigation     & 15     & 135.17 & 222.52   & 163.05 & 224.64  & 125.26 & \textbf{212.04} & 130.36 & 216.18 \\
                   & 300    & 151.27 & 1015.66  & 371.62 & 978.06  & 125.26 & 955.57          & 151.98 & \textbf{820.67} \\
    Many Obstacles & 15     & 243.10 & 424.24   & 596.45 & 446.54  & 239.39 & \textbf{417.21} & 250.69 & 430.88 \\
                   & 3000   & 265.06 & 13708.01 & 1471.25 & 8120.80 & 239.39 & 14811.65        & 456.85 & \textbf{7579.74} \\
    Pendulum       & 1      & 554.05 & 639.20   & 688.35 & 679.30  & 548.53 & 638.33          & 558.30 & \textbf{606.81} \\
                   & 10     & 756.96 & 817.84   & 743.83 & 738.78  & 548.53 & 722.23          & 562.10 & \textbf{709.60} \\
    Planar Arm     & 3      & 273.10 & 285.36   & 184.77 & \textbf{205.71} & 189.96 & 217.43 & 204.63 & 227.20 \\
                   & 300    & 694.06 & 1123.56  & 332.81 & 789.72  & 189.96 & 837.36          & 288.58 & \textbf{762.19} \\
    \bottomrule
  \end{tabular}
\end{table}

\subsection{Parameter and Horizon Sensitivity}
\label{subsec:parameter_sensitivity}

We assess whether CBO-MPPI's free-energy advantage over MPPI persists beyond the reference settings, varying first the sampling parameters and then the prediction horizon. We vary the inverse temperature $\alpha$ and the ratio $\kappa=D\sigma_{\rm cbo}^2/(2\lambda)$ one at a time, with 20 paired trials per setting. Across 80 task--temperature--parameter settings, the median paired $F_{\rm CL}$ difference (CBO-MPPI minus MPPI) is negative at 62, including all 40 high-temperature settings. Thus the high-temperature advantage is not tied to the reference parameter choice. At low temperature, several settings reverse the comparison, and Pendulum at $\tau=1$ favors MPPI throughout.

The horizon sweep shows where consensus helps. At short horizons, MS-MPPI can match or outperform CBO-MPPI. At the longest tested horizons, however, CBO-MPPI succeeds in $17/20$ trials versus MS-MPPI's $11/20$ on Navigation and $16/20$ versus MS-MPPI's $1/20$ on Pendulum. At the reference horizons, CBO-MPPI has lower realized cost and higher success rates than MS-MPPI on all three tasks; its free energy is lower on Navigation and Planar Arm and nearly equal on Pendulum. \fullversionref{Appendices~\ref{full:subsec:horizon_effect} and~\ref{full:subsec:parameter_sensitivity} give}{The supplementary material gives} the grids and complete results.

\FloatBarrier

\section{Conclusion}
\label{sec:conclusion}

We proposed CBO-MPPI for predictive-control problems in which hard constraints can make the trajectory cost discontinuous and create multiple search basins. Gaussian perturbations transform this cost into an MPPI free energy that rewards sampling neighborhoods containing feasible, low-cost trajectories. CBO-MPPI searches this smoothed objective with a population of coupled centers and obtains the control from the selected center's MPPI weighted mean. For the mean-field model of a fixed open-loop problem with exact, unprojected free energy, we established well-posedness and, under additional assumptions, convergence toward a free-energy minimizer. The simulations support this design: when the temperature is high enough for perturbation feasibility to matter, the consensus search reaches lower free energy than MPPI and MS-MPPI and keeps trajectories farther from obstacles than CBO applied to the nominal cost, and its advantage over MPPI persists across the tested parameter ranges. Extending the theory to finite-particle and finite-rollout implementations, input projection, and closed-loop stability remains future work.

%% file: endmatter.tex
\FloatBarrier
\subsection*{Generative AI Use}
Generative AI tools were used to assist with numerical computations and writing. The authors reviewed and verified all AI-assisted content and take responsibility for the final content of the paper.
\acks{D. Kalise is supported by the Air Force Office of Scientific Research under award number FA8655-26-1-B011. Y. Huang is supported by a Research Fellowship from the Royal Commission for the Exhibition of 1851.}
\bibliography{references}

%% file: appendix.tex
\section{Mean-Field Convergence and Proofs}
\label{sec:proof}

This appendix proceeds in two steps. Subsection~\ref{subsec:proof_free_energy_properties} establishes the free-energy properties needed for well-posedness and proves Theorem~\ref{prop:mean_field_well_posedness}. Subsection~\ref{subsec:paper_global_convergence} introduces the additional assumptions and objective conditions under which Theorem~\ref{thm:mean_field_global_convergence} establishes finite-accuracy global convergence.

\subsection{Free-Energy Properties and Well-Posedness}
\label{subsec:proof_free_energy_properties}
\label{subsec:proof_mean_field_well_posedness}

We first verify the regularity and growth properties of $F_x$ needed for well-posedness and then prove Theorem~\ref{prop:mean_field_well_posedness}. Additional conditions for finite-accuracy global convergence are introduced in Subsection~\ref{subsec:paper_global_convergence}.
We begin by introducing the partition function
\[
  Z(M)
  :=
  \int_{\mathcal C}
  \varphi_\Sigma(U-M)
  \exp\left(
  -\frac{J_x(U)}{\tau}
  \right)dU,
\]
so that $F_x(M)=-\tau\log Z(M)$.

\begin{lemma}[Positivity and finiteness of the partition function]
  \label{lem:partition_positive_finite}
  Under Assumptions~\ref{ass:admissible_set} and \ref{ass:cost}, for every $M\in\mathbb{R}^D$,
  \[
    0<Z(M)<\infty.
  \]
  Consequently, $F_x(M)$ is well-defined for all $M\in\mathbb{R}^D$.
\end{lemma}

\begin{proof}
  By Assumption~\ref{ass:cost},
  \[
    e^{-J_{\max}/\tau}
    \leq
    \exp\left(-\frac{J_x(U)}{\tau}\right)
    \leq
    e^{-J_{\min}/\tau},
    \qquad U\in\mathcal C.
  \]
  Since the Gaussian density $\varphi_\Sigma$ is everywhere positive and $\mathcal C$ has positive Lebesgue measure,
  \[
    Z(M)
    \geq
    e^{-J_{\max}/\tau}
    \int_{\mathcal C}
    \varphi_\Sigma(U-M)dU
    >
    0.
  \]
  On the other hand,
  \[
    Z(M)
    \leq
    e^{-J_{\min}/\tau}
    \int_{\mathcal C}
    \varphi_\Sigma(U-M)dU
    \leq
    e^{-J_{\min}/\tau}
    <
    \infty,
  \]
  which completes the proof.
\end{proof}
We next show that Gaussian convolution makes $F_x$ smooth, even when $J_x$ is discontinuous.

\begin{lemma}[Smoothness of the MPPI free energy]
  \label{lem:free_energy_smoothness}
  Under Assumptions~\ref{ass:admissible_set} and \ref{ass:cost},
  \[
    Z\in C^\infty(\mathbb{R}^D),
    \qquad
    F_x\in C^\infty(\mathbb{R}^D).
  \]
\end{lemma}

\begin{proof}
  Set
  \[
    h(U)
    :=
    \mathbf{1}_{\mathcal C}(U)
    \exp\left(-\frac{J_x(U)}{\tau}\right).
  \]
  By Assumption~\ref{ass:cost}, $h$ is a bounded measurable function, and by Assumption~\ref{ass:admissible_set} it has compact support.
  Then
  \[
    Z(M)
    =
    \int_{\mathbb{R}^D}
    \varphi_\Sigma(U-M)h(U)dU.
  \]
  The Gaussian density $\varphi_\Sigma$ belongs to $C^\infty$, and each of its derivatives is a product of a polynomial and a Gaussian function.
  Since $h$ is bounded with compact support, differentiation and integration may be interchanged for any multi-index $\beta$:
  \[
    \partial_M^\beta Z(M)
    =
    \int_{\mathbb{R}^D}
    \partial_M^\beta
    \varphi_\Sigma(U-M)
    h(U)dU.
  \]
  Hence $Z\in C^\infty(\mathbb{R}^D)$.
  Since $Z(M)>0$ by Lemma~\ref{lem:partition_positive_finite}, the function $\log Z(M)$ is also $C^\infty$, and therefore
  \[
    F_x(M)=-\tau\log Z(M)\in C^\infty(\mathbb{R}^D).
  \]
\end{proof}

For $M\in\mathbb R^D$, define the tilted distribution
\[
  \rho_M(dU)
  :=
  \frac{
    \varphi_\Sigma(U-M)
    \exp\left(-J_x(U)/\tau\right)
    \mathbf 1_{\mathcal C}(U)\,dU
  }{Z(M)}.
\]
Then $\mathbb E_{\rho_M}[U]=\Phi(M)$, where $\Phi$ is the exact MPPI update already defined in~\eqref{eq:mppi_gradient_identity}.
We recall the Gaussian gradient and Hessian identities of \cite[Lemma~1]{Fazlyab2026Model}, whose $\theta=\mu$, $u$, $f_0$, $C$, and $\pi_\theta(u)$ correspond to our $M$, $U$, $J_x$, $\mathcal C$, and $\varphi_\Sigma(U-M)$, respectively.
That reference assumes a continuous cost, but its derivation uses only positivity of $Z$ and differentiation under the integral, which Lemmas~\ref{lem:partition_positive_finite} and~\ref{lem:free_energy_smoothness} provide for bounded Borel $J_x$.

\begin{proposition}[Gaussian free-energy identities \eqref{eq:mppi_gradient_identity}]
  \label{prop:mppi_preconditioned_gradient}
  Under Assumptions~\ref{ass:admissible_set} and \ref{ass:cost},
  \begin{equation}
    \nabla F_x(M)
    =
    \tau\Sigma^{-1}
    \left(M-\mathbb E_{\rho_M}[U]\right),
    \label{eq:free_energy_gradient}
  \end{equation}
  and
  \begin{equation}
    \nabla^2F_x(M)
    =
    \tau\Sigma^{-1}
    -
    \tau\Sigma^{-1}
    \operatorname{Cov}_{\rho_M}(U)
    \Sigma^{-1}.
    \label{eq:free_energy_hessian_covariance}
  \end{equation}
\end{proposition}

\begin{proof}
  See \cite[Lemma~1]{Fazlyab2026Model}.
\end{proof}

We next show that $F_x$ is coercive, i.e., $F_x(M)\to\infty$ as $\|M\|\to\infty$. This ensures that the objective does not attain values near its global minimum far from the minimizer and implies the existence of a global minimizer.
We use the $\Sigma^{-1}$-weighted distance notation
\[
  \|z\|_{\Sigma^{-1}}^2:=z^\top\Sigma^{-1}z,
  \qquad \operatorname{dist}_{\Sigma^{-1}}(M,\mathcal C)
  :=\inf_{U\in\mathcal C}\|M-U\|_{\Sigma^{-1}}.
\]

\begin{lemma}[Coercivity of the MPPI free energy]
  \label{lem:free_energy_coercive}
  Under Assumptions~\ref{ass:admissible_set} and \ref{ass:cost}, there exists a constant $c_0\in\mathbb{R}$ such that for all $M\in\mathbb{R}^D$,
  \begin{equation}
    F_x(M)
    \geq
    \frac{\tau}{2}
    \operatorname{dist}_{\Sigma^{-1}}(M,\mathcal C)^2
    -
    c_0.
    \label{eq:free_energy_coercive_bound}
  \end{equation}
  In particular,
  \[
    F_x(M)\to+\infty
    \qquad
    \text{as } \|M\|\to\infty.
  \]
\end{lemma}

\begin{proof}
  The Gaussian density satisfies
  \[
    \varphi_\Sigma(U-M)
    =\frac{1}{(2\pi)^{D/2}\det(\Sigma)^{1/2}}
    \exp\left(-\frac12\|U-M\|_{\Sigma^{-1}}^2\right).
  \]
  For any $U\in\mathcal C$,
  \[
    \|U-M\|_{\Sigma^{-1}}^2
    \geq
    \operatorname{dist}_{\Sigma^{-1}}(M,\mathcal C)^2,
  \]
  so
  \[
    \varphi_\Sigma(U-M)
    \leq\frac{1}{(2\pi)^{D/2}\det(\Sigma)^{1/2}}
    \exp\left(-\frac12\operatorname{dist}_{\Sigma^{-1}}(M,\mathcal C)^2\right).
  \]
  Combined with the bound $\exp(-J_x(U)/\tau)\leq\exp(-J_{\min}/\tau)$ from Assumption~\ref{ass:cost},
  \[
    \begin{aligned}
      Z(M)
       & =
      \int_{\mathcal C}
      \varphi_\Sigma(U-M)
      \exp\left(-\frac{J_x(U)}{\tau}\right)dU
      \\
       & \leq
      \frac{|\mathcal C|}{(2\pi)^{D/2}\det(\Sigma)^{1/2}}
      \exp\left(-\frac{J_{\min}}{\tau}\right)
      \exp\left(
      -\frac12
      \operatorname{dist}_{\Sigma^{-1}}(M,\mathcal C)^2
      \right),
    \end{aligned}
  \]
  where $|\mathcal C|$ denotes the Lebesgue measure of $\mathcal C$.
  Applying $-\tau\log$ to both sides,
  \[
    F_x(M)
    =
    -\tau\log Z(M)
    \geq
    \frac{\tau}{2}
    \operatorname{dist}_{\Sigma^{-1}}(M,\mathcal C)^2
    -
    c_0,
  \]
  where $c_0$ depends on $\mathcal C$, $\Sigma$, $\tau$, and $J_{\min}$.
  Since $\mathcal C$ is compact,
  \[
    \|M\|\to\infty
    \quad\Longrightarrow\quad
    \operatorname{dist}_{\Sigma^{-1}}(M,\mathcal C)\to\infty,
  \]
  and therefore $F_x(M)\to+\infty$.
\end{proof}

\begin{corollary}[Existence of a global minimizer]
  \label{cor:existence_global_minimizer}
  Under Assumptions~\ref{ass:admissible_set} and \ref{ass:cost}, $F_x$ has at least one global minimizer; that is,
  \[
    \operatorname*{argmin}_{M\in\mathbb{R}^D}F_x(M)
    \neq
    \emptyset.
  \]
\end{corollary}

\begin{proof}
  $F_x$ is continuous by Lemma~\ref{lem:free_energy_smoothness} and coercive by Lemma~\ref{lem:free_energy_coercive}.
  Hence there exists $R>0$ large enough that
  \[
    \inf_{M\in\mathbb R^D}F_x(M)
    =
    \inf_{\|M\|\leq R}F_x(M).
  \]
  Since the closed ball $\{M\in\mathbb R^D:\|M\|\leq R\}$ is compact, the Weierstrass theorem guarantees that the continuous function $F_x$ attains its minimum on this ball.
\end{proof}

Fix a minimizer
\[
  M^\star\in\operatorname*{argmin}_{M\in\mathbb R^D}F_x(M).
\]
The Gaussian free-energy structure also provides the global regularity and growth estimates required by the standard CBO well-posedness theory.

\begin{lemma}[Global regularity and quadratic growth of the MPPI free energy]
  \label{lem:free_energy_global_regularity}
  Under Assumptions~\ref{ass:admissible_set} and \ref{ass:cost}, there exists $L_F<\infty$ such that
  \[
    \sup_{M\in\mathbb R^D}\|\nabla^2F_x(M)\|\leq L_F.
  \]
  For this fixed global minimizer $M^\star$, define
  \[
    \widetilde F_x(z):=F_x(M^\star+z)-F_x(M^\star).
  \]
  Then, for all $z,y\in\mathbb R^D$,
  \begin{equation}
    0\leq\widetilde F_x(z)\leq \frac{L_F}{2}\|z\|^2,
    \qquad |\widetilde F_x(z)-\widetilde F_x(y)|
    \leq L_F(\|z\|+\|y\|)\|z-y\|.
    \label{eq:translated_free_energy_upper_conditions}
  \end{equation}
  Moreover, there exist $R>0$ and $c_l>0$ such that
  \begin{equation}
    \widetilde F_x(z)\geq c_l\|z\|^2,
    \qquad \|z\|\geq R.
    \label{eq:translated_free_energy_lower_growth}
  \end{equation}
\end{lemma}

\begin{proof}
  We use the Hessian identity recalled in Proposition~\ref{prop:mppi_preconditioned_gradient}. Since $\rho_M$ is supported on the compact set $\mathcal C$, choose $R_{\mathcal C,0}<\infty$ such that $\|U\|\leq R_{\mathcal C,0}$ on $\mathcal C$. For every unit vector $v$,
  \[
    v^\top\operatorname{Cov}_{\rho_M}(U)v
    =\operatorname{Var}_{\rho_M}(v^\top U)
    \leq\mathbb E_{\rho_M}|v^\top U|^2
    \leq R_{\mathcal C,0}^2.
  \]
  Thus $\|\operatorname{Cov}_{\rho_M}(U)\|\leq R_{\mathcal C,0}^2$ uniformly in $M$, and \eqref{eq:free_energy_hessian_covariance} gives, for example,
  \[
    L_F
    :=
    \tau\|\Sigma^{-1}\|
    +\tau R_{\mathcal C,0}^2\|\Sigma^{-1}\|^2.
  \]

  Since $F_x$ is smooth and $M^\star$ is a minimizer, $\widetilde F_x(0)=0$ and $\nabla\widetilde F_x(0)=0$. The global Hessian bound therefore yields
  \[
    \|\nabla\widetilde F_x(z)\|\leq L_F\|z\|.
  \]
  Taylor's formula at the origin gives the quadratic upper bound in \eqref{eq:translated_free_energy_upper_conditions}. Integrating the gradient along the segment from $y$ to $z$ gives
  \[
    |\widetilde F_x(z)-\widetilde F_x(y)|
    \leq L_F\int_0^1\|y+t(z-y)\|dt\,\|z-y\|
    \leq L_F(\|z\|+\|y\|)\|z-y\|,
  \]
  which is the stated weighted-Lipschitz estimate.

  Finally, let
  \[
    R_{\mathcal C}:=\sup_{U\in\mathcal C}\|U-M^\star\|<\infty.
  \]
  Whenever $\|z\|\geq R_{\mathcal C}$, the reverse triangle inequality gives, for every $U\in\mathcal C$,
  \[
    \|M^\star+z-U\|_{\Sigma^{-1}}^2
    \geq
    \lambda_{\min}(\Sigma^{-1})(\|z\|-R_{\mathcal C})^2.
  \]
  Combining this inequality with \eqref{eq:free_energy_coercive_bound} gives
  \[
    \widetilde F_x(z)
    \geq
    \frac{\tau}{2}\lambda_{\min}(\Sigma^{-1})
    (\|z\|-R_{\mathcal C})^2-c_0-F_x(M^\star).
  \]
  For sufficiently large $\|z\|$, the constants can be absorbed into the quadratic term, yielding \eqref{eq:translated_free_energy_lower_growth}.
\end{proof}

\paragraph{Proof of well-posedness.}

The law $\rho_s:=\operatorname{Law}(M_s)$ of the process in Theorem~\ref{prop:mean_field_well_posedness} is a weak solution of the nonlinear Fokker--Planck equation
\begin{equation}
  \partial_s\rho_s
  =\lambda\nabla_M\cdot\left[(M-m_\alpha[\rho_s])\rho_s\right]
  +\frac{\sigma_{\mathrm{cbo}}^2}{2}
  \Delta_M\left[\|M-m_\alpha[\rho_s]\|^2\rho_s\right].
  \label{eq:mean_field_fokker_planck}
\end{equation}

For completeness, the weak formulation of \eqref{eq:mean_field_fokker_planck} is as follows: for every $\phi\in C_c^\infty(\mathbb R^D)$ and $s\in[0,S]$,
\begin{equation}
  \begin{aligned}
    &\int_{\mathbb R^D}\phi(M)\rho_s(dM)
    -\int_{\mathbb R^D}\phi(M)\rho_0(dM)\\
    &=-\lambda\int_0^s\int_{\mathbb R^D}
      (M-m_\alpha[\rho_{s'}])\cdot\nabla\phi(M)\rho_{s'}(dM)\,ds'\\
    &\quad+\frac{\sigma_{\mathrm{cbo}}^2}{2}
      \int_0^s\int_{\mathbb R^D}
      \|M-m_\alpha[\rho_{s'}]\|^2\Delta\phi(M)\rho_{s'}(dM)\,ds'.
  \end{aligned}
  \label{eq:mean_field_fokker_planck_weak}
\end{equation}

\begin{proof}[Proof of Theorem~\ref{prop:mean_field_well_posedness}]
  Use the global minimizer $M^\star$ fixed after Corollary~\ref{cor:existence_global_minimizer}, translate $z=M-M^\star$, and subtract the constant $F_x(M^\star)$ from the objective. Lemma~\ref{lem:free_energy_global_regularity} verifies the boundedness-from-below, weighted-Lipschitz, and quadratic upper-growth conditions in \cite[Assumption~3.1]{Carrillo2016analytical}, as well as the far-field quadratic lower-growth condition used in \cite[Theorem~3.2]{Carrillo2016analytical}. The coordinate translation only translates the consensus point, and the additive objective constant cancels from its normalized Boltzmann ratio. Applying \cite[Theorem~3.2]{Carrillo2016analytical} therefore yields a unique nonlinear strong process on every finite interval, whose law satisfies \eqref{eq:mean_field_fokker_planck_weak}.

  The proof of \cite[Theorem~3.2]{Carrillo2016analytical} retains the fourth-moment estimate established in Step~2 of the proof of \cite[Theorem~3.1]{Carrillo2016analytical}, which gives the fourth-moment bound in Theorem~\ref{prop:mean_field_well_posedness}; it also bounds $m_\alpha[\rho_s]$ uniformly on $[0,S]$. Consequently,
  \begin{equation}
    \sup_{0\leq s\leq S}
    \mathbb E\|M_s-m_\alpha[\rho_s]\|^4
    <\infty.
    \label{eq:mean_field_centered_fourth_moment}
  \end{equation}
  For $0\leq s<t\leq S$, the same process provides a coupling of $\rho_s$ and $\rho_t$. Applying Jensen's inequality to the drift integral and the Burkholder--Davis--Gundy inequality to the stochastic integral in \eqref{eq:mean_field_cbo_mppi}, together with \eqref{eq:mean_field_centered_fourth_moment}, gives
  \[
    W_4^4(\rho_t,\rho_s)
    \leq
    \mathbb E\|M_t-M_s\|^4
    \leq
    C_S\left(|t-s|^4+|t-s|^2\right).
  \]
  Thus $W_4(\rho_t,\rho_s)\to0$ as $t\to s$, strengthening the $\mathcal P_2$ continuity stated in \cite[Theorem~3.2]{Carrillo2016analytical} to $\rho\in C([0,S],\mathcal P_4(\mathbb R^D))$.
\end{proof}

\subsection{Additional Assumptions and Finite-Accuracy Global Convergence}
\label{subsec:paper_global_convergence}

We apply the isotropic CBO convergence theorem \cite[Theorem~3.7]{Fornasier2024ConsensusBased} to the mean-field model in Section~\ref{sec:theoretical_analysis}.
We retain the fixed $x,H,\tau>0,\Sigma\succ0$ and exact, unprojected free energy $F_x$ used there.
We write $W_2(\mu,\nu)$ for the 2-Wasserstein distance between probability measures $\mu$ and $\nu$, and $\delta_M$ for the Dirac measure at $M$.
The existence of a global minimizer follows from Corollary~\ref{cor:existence_global_minimizer}. For finite-accuracy global convergence, we additionally assume uniqueness and nondegeneracy.

\begin{assumption}[Unique nondegenerate global minimizer]
  \label{ass:unique_nondegenerate_minimizer}
  The global minimizer $M^\star$ of $F_x$ is unique and satisfies $\nabla^2F_x(M^\star)\succ0$.
\end{assumption}
Assumption~\ref{ass:unique_nondegenerate_minimizer} is not automatic for general MPC problems; symmetric obstacle avoidance, for example, may yield multiple equivalent minimizers. A different analysis may permit a set-valued conclusion: for a modified CBO model, \citet{Fornasier2025PDE} prove exponential concentration toward a set of global minimizers with multiple smooth, convex, and compact components. We assume uniqueness to apply the standard theorem directly.
Here and throughout, $B_r(M^\star):=\{M\in\mathbb R^D:\|M-M^\star\|<r\}$.

The next lemma is the standard Taylor estimate at a nondegenerate minimizer. We apply it under Assumption~\ref{ass:unique_nondegenerate_minimizer} to verify the inverse-continuity condition in \cite[Definition~3.5]{Fornasier2024ConsensusBased}.

\begin{lemma}[Local quadratic growth near the nondegenerate minimizer]
  \label{lem:local_quadratic_growth}
  Suppose Assumptions~\ref{ass:admissible_set}, \ref{ass:cost}, and \ref{ass:unique_nondegenerate_minimizer} hold.
  Then there exist constants $r>0$ and $\mu>0$ such that for all $M\in B_r(M^\star)$,
  \begin{equation}
    F_x(M)-F_x(M^\star)
    \geq
    \frac{\mu}{2}\|M-M^\star\|^2.
    \label{eq:local_quadratic_growth}
  \end{equation}
  Consequently,
  \begin{equation}
    \|M-M^\star\|
    \leq
    \sqrt{\frac{2}{\mu}}
    \left(
    F_x(M)-F_x(M^\star)
    \right)^{1/2}
    \label{eq:inverse_continuity}
  \end{equation}
  holds on $B_r(M^\star)$.
\end{lemma}

\begin{proof}
  By Assumption~\ref{ass:unique_nondegenerate_minimizer}, $\nabla^2F_x(M^\star)\succ0$, so there exists $\mu>0$ such that $\nabla^2F_x(M^\star)\succeq 2\mu I_D$.
  Since $F_x\in C^\infty$ by Lemma~\ref{lem:free_energy_smoothness}, $\nabla^2F_x$ is continuous, and for $r>0$ sufficiently small, $\nabla^2F_x(M)\succeq \mu I_D$ for all $M\in B_r(M^\star)$.
  Since $F_x$ is smooth and $M^\star$ is a global minimizer, $\nabla F_x(M^\star)=0$.
  Taylor's theorem then gives, for any $M\in B_r(M^\star)$, some $\theta\in(0,1)$, and $M_\theta:=M^\star+\theta(M-M^\star)$,
  \[
    \begin{aligned}
      F_x(M)-F_x(M^\star)
       & =
      \nabla F_x(M^\star)\cdot(M-M^\star)
      +\frac12
      (M-M^\star)^\top
      \nabla^2F_x(M_\theta)
      (M-M^\star)
      \\
       & \geq
      \frac{\mu}{2}\|M-M^\star\|^2,
    \end{aligned}
  \]
  establishing \eqref{eq:local_quadratic_growth}. Rearranging gives \eqref{eq:inverse_continuity}.
\end{proof}

\begin{lemma}[CBO objective conditions]
  \label{lem:cbo_objective_conditions}
  Under Assumptions~\ref{ass:admissible_set}, \ref{ass:cost}, and \ref{ass:unique_nondegenerate_minimizer}, there exist constants $c_{\mathrm{ic}},R_0,E_\infty>0$ such that
  \begin{equation}
    \|M-M^\star\|
    \leq\frac{1}{c_{\mathrm{ic}}}
    \left(F_x(M)-F_x(M^\star)\right)^{1/2},
    \qquad M\in B_{R_0}(M^\star),
    \label{eq:cbo_inverse_continuity_global_ball}
  \end{equation}
  and
  \begin{equation}
    F_x(M)-F_x(M^\star)>E_\infty,
    \qquad
    M\in B_{R_0}(M^\star)^c.
    \label{eq:cbo_far_field_gap}
  \end{equation}
\end{lemma}

\begin{proof}
  By Lemma~\ref{lem:local_quadratic_growth}, there exist $r>0$ and $\mu>0$ such that on $B_r(M^\star)$,
  \[
    \|M-M^\star\|
    \leq
    \sqrt{\frac{2}{\mu}}
    \left(F_x(M)-F_x(M^\star)\right)^{1/2}.
  \]
  By Lemma~\ref{lem:free_energy_coercive}, $F_x(M)\to\infty$ as $\|M\|\to\infty$.
  Choosing $R_0>r$ large enough gives
  \[
    \inf_{M\in B_{R_0}(M^\star)^c}
    \left(F_x(M)-F_x(M^\star)\right)
    >
    0.
  \]
  Setting
  \[
    E_\infty
    :=
    \frac12
    \inf_{M\in B_{R_0}(M^\star)^c}
    \left(F_x(M)-F_x(M^\star)\right)
    >0
  \]
  yields \eqref{eq:cbo_far_field_gap}.
  For the annular region
  \[
    A
    :=
    \{M\in\mathbb R^D:
    r\leq \|M-M^\star\|\leq R_0\},
  \]
  $A$ is compact, $F_x$ is continuous, and $M^\star\notin A$, so
  \[
    \delta_A
    :=
    \min_{M\in A}
    \left(F_x(M)-F_x(M^\star)\right)
    >
    0.
  \]
  For $M\in A$,
  \[
    \|M-M^\star\|
    \leq
    R_0
    \leq
    \frac{R_0}{\sqrt{\delta_A}}
    \left(F_x(M)-F_x(M^\star)\right)^{1/2}.
  \]
  Setting
  \[
    c_{\mathrm{ic}}
    :=
    \min
    \left\{
    \sqrt{\frac{\mu}{2}},
    \frac{\sqrt{\delta_A}}{R_0}
    \right\}
  \]
  gives \eqref{eq:cbo_inverse_continuity_global_ball} uniformly on $B_{R_0}(M^\star)$.
\end{proof}

Global convergence also requires an initial support condition and a bound on the diffusion strength.

\begin{assumption}[Initial distribution and CBO parameters]
  \label{ass:cbo_initial_and_parameters}
  The initial distribution $\rho_0\in\mathcal P_4(\mathbb{R}^D)$ satisfies $\rho_0(B_r(M^\star))>0$ for every $r>0$.
  Furthermore, $\lambda,\sigma_{\mathrm{cbo}}>0$ satisfy $2\lambda-D\sigma_{\mathrm{cbo}}^2>0$.
\end{assumption}
The condition $\rho_0(B_r(M^\star))>0$ for every $r>0$ is equivalent to $M^\star\in\operatorname{supp}\rho_0$. That is, $\rho_0$ assigns positive probability to every neighborhood of the global minimizer.

We analyze finite-accuracy convergence using the energy functional
\[
  \mathcal V(\rho)
  :=
  \frac12\int_{\mathbb R^D}\|M-M^\star\|^2\rho(dM),
\]
which measures the mean squared distance from the unique global minimizer $M^\star$.

\begin{theorem}[Finite-accuracy global convergence; \citealp{Fornasier2024ConsensusBased}, Theorem~3.7]
  \label{thm:mean_field_global_convergence}
  Suppose Assumptions~\ref{ass:admissible_set}, \ref{ass:cost},
  \ref{ass:unique_nondegenerate_minimizer},
  and \ref{ass:cbo_initial_and_parameters} hold.
  For any
  $\varepsilon\in(0,\mathcal V(\rho_0))$ and any
  $\vartheta\in(0,1)$,
  define
  \begin{equation}
    S^\star
    :=
    \frac{1}{(1-\vartheta)(2\lambda-D\sigma_{\mathrm{cbo}}^2)}
    \log\left(\frac{\mathcal V(\rho_0)}{\varepsilon}\right).
    \label{eq:cbo_accuracy_time}
  \end{equation}
  There exists $\alpha_0>0$, depending in particular on $\varepsilon$, $\vartheta$, $D$, $\lambda$, $\sigma_{\mathrm{cbo}}$, the initial law $\rho_0$, and the objective constants and local sublevel geometry of $F_x$, such that for every $\alpha>\alpha_0$, the law $\rho_s:=\operatorname{Law}(M_s)$ of the unique nonlinear process $M$ from Theorem~\ref{prop:mean_field_well_posedness}, initialized with law $\rho_0$, satisfies
  \[
    \min_{0\le s\le S^\star}\mathcal V(\rho_s)\leq\varepsilon.
  \]
  Moreover, until $\mathcal V(\rho_s)$ first reaches $\varepsilon$,
  \begin{equation}
    \mathcal V(\rho_s)
    \leq
    \mathcal V(\rho_0)
    \exp\left(-(1-\vartheta)(2\lambda-D\sigma_{\mathrm{cbo}}^2)s\right).
    \label{eq:cbo_exponential_decay}
  \end{equation}
  In particular, there exists $s_\varepsilon\in[0,S^\star]$ such that $W_2(\rho_{s_\varepsilon},\delta_{M^\star})\leq\sqrt{2\varepsilon}$.
\end{theorem}

\begin{proof}
  Lemmas~\ref{lem:partition_positive_finite} and
  \ref{lem:free_energy_smoothness} show that $F_x$ is a continuous real-valued
  objective on $\mathbb R^D$. Lemma~\ref{lem:free_energy_coercive} and
  Corollary~\ref{cor:existence_global_minimizer} show that it attains its
  global minimum. Hence condition A1 in
  \cite[Definition~3.5]{Fornasier2024ConsensusBased} holds.

  Under Assumption~\ref{ass:unique_nondegenerate_minimizer},
  Lemma~\ref{lem:cbo_objective_conditions} gives
  \eqref{eq:cbo_inverse_continuity_global_ball} and
  \eqref{eq:cbo_far_field_gap}. These are the inverse-continuity and
  far-field parts of condition A2 in
  \cite[Definition~3.5]{Fornasier2024ConsensusBased}, with exponent
  $\nu=1/2$.

  Assumption~\ref{ass:cbo_initial_and_parameters} gives
  $\rho_0\in\mathcal P_4(\mathbb R^D)$,
  $M^\star\in\operatorname{supp}\rho_0$, and
  $2\lambda>D\sigma_{\mathrm{cbo}}^2$. Theorem~\ref{prop:mean_field_well_posedness} supplies the unique process and its law $\rho\in C([0,S^\star],\mathcal P_4(\mathbb R^D))$. Equation~\eqref{eq:mean_field_fokker_planck_weak} shows that this law weakly solves \eqref{eq:mean_field_fokker_planck}.

  The dynamics \eqref{eq:mean_field_cbo_mppi} coincide with the standard isotropic CBO dynamics of \cite[Theorem~3.7]{Fornasier2024ConsensusBased} when the cutoff function is identically equal to one. The probability-mass estimate \cite[Proposition~4.6]{Fornasier2024ConsensusBased}, together with $M^\star\in\operatorname{supp}\rho_0$ and finite-time boundedness of $m_\alpha[\rho_s]$, gives, for each fixed $r>0$ and $S>0$,
  \[
  \rho_s(B_r(M^\star))
  \geq\left(\int_{\mathbb R^D}\phi_r(M)\rho_0(dM)\right)
  \times e^{-p_{r,S}s},\qquad 0\leq s\leq S,
  \]
  where $\phi_r$ is the compactly supported mollifier in the cited proposition and $p_{r,S}$ is a finite constant whose dependence includes the finite-time consensus bound. Since $\phi_r$ is bounded below by some $c_r>0$ on $B_{r/2}(M^\star)$,
  \[
    \int_{\mathbb R^D}\phi_r(M)\rho_0(dM)
    \geq
    c_r\rho_0(B_{r/2}(M^\star))
    >0.
  \]
  Hence $\inf_{0\leq s\leq S}\rho_s(B_r(M^\star))>0$. This lower bound controls the denominator in the quantitative Laplace estimate underlying \cite[Theorem~3.7]{Fornasier2024ConsensusBased}.

  Applying \cite[Theorem~3.7]{Fornasier2024ConsensusBased} with dimension $d=D$, objective $\mathcal E=F_x$, and diffusion coefficient $\sigma=\sigma_{\mathrm{cbo}}$ gives the horizon
  \eqref{eq:cbo_accuracy_time}, the exponential estimate
  \eqref{eq:cbo_exponential_decay}, and a time
  $s_\varepsilon\leq S^\star$ at which
  $\mathcal V(\rho_{s_\varepsilon})\leq\varepsilon$.
  Since
  $W_2^2(\rho_s,\delta_{M^\star})=2\mathcal V(\rho_s)$, the final bound
  follows.
\end{proof}

The theorem establishes concentration near the global minimizer of $F_x$. Under the stated assumptions, any prescribed tolerance $\varepsilon$ can be reached by choosing $\alpha$ sufficiently large. At some CBO pseudo-time $s_\varepsilon\leq S^\star$, the root mean squared distance from $M^\star$ is at most $\sqrt{2\varepsilon}$, with an exponential bound on the mean squared distance up to this first attainment. The result applies to a fixed open-loop problem with exact, unprojected free energy and mean-field dynamics. Finite-particle and finite-rollout convergence, input projection, warm starts, and closed-loop stability remain open.

\paperfloatbarrier

\section{Algorithms and Numerical Details}
\label{sec:paper_full_details}

This appendix documents the three comparison controllers. Subsection~\ref{full:sec:proposed_method} specifies MPPI, MS-MPPI, and CBO-MPC, including their output and warm-start rules.

\subsection{Baseline Controller Algorithms}
\label{full:sec:proposed_method}

\subsubsection{MPPI}
\label{subsec:mppi_free_energy}
MPPI maintains one nominal sequence and applies the weighted-mean update in~\eqref{eq:sample_mppi_update}. The updated sequence supplies the first input and is shifted for the next control period.

\begin{algorithm}[htbp]
\algocfconts{alg:mppi}{\caption{MPPI: one control period}}{
\begin{algorithmic}[1]
\REQUIRE State $x$; initialized or warm-started sequence $M\in\R^{mH}$; $L$, $\tau>0$, $\Sigma\succ0$
\ENSURE Applied input and next-period nominal sequence
\STATE Draw $\epsilon_j\sim\mathcal N(0,\Sigma)$ and set $U_j=M+\epsilon_j$, $j=1,\ldots,L$
\STATE Roll out from $x$ and evaluate $S_j=\bar J_x(U_j)$ for all $j$ in parallel
\STATE Compute $\widehat Z=L^{-1}\sum_j e^{-S_j/\tau}$ and $M^+$ by~\eqref{eq:sample_mppi_update}, using the zero-weight safeguard if needed
\STATE Apply $[M^+]_0$ and set $M^{\mathrm{next}}=\operatorname{Shift}(M^+)$
\end{algorithmic}}
\end{algorithm}
\FloatBarrier

\subsubsection{MS-MPPI}
\label{subsec:paper_multistart}
MS-MPPI maintains $N$ persistent centers $M_i^k\in\R^D$. At each round, it draws $U_{i,j}^k=M_i^k+\epsilon_{i,j}^k$, with $\epsilon_{i,j}^k\sim\mathcal N(0,\Sigma)$, and evaluates $S_{i,j}^k=\bar J_x(U_{i,j}^k)$. Each local batch gives
\begin{equation}
\widehat Z_i^k=\frac1L\sum_{j=1}^L e^{-S_{i,j}^k/\tau},\qquad
\widehat F_i^k=-\tau\log\widehat Z_i^k,\qquad
w_{i,j}^k=\frac{e^{-S_{i,j}^k/\tau}}{L\widehat Z_i^k}.
\label{eq:paper_multistart_weights}
\end{equation}
For the first $K-1$ rounds, every center follows its own MPPI update,
\begin{equation}
M_i^{k+1}=\sum_{j=1}^L w_{i,j}^k U_{i,j}^k,\qquad k=0,\ldots,K-2.
\label{eq:paper_multistart_update}
\end{equation}
The final round selects $i^\star\in\argmin_i\widehat F_i^{K-1}$ and returns $U_{\rm out}=\sum_j w_{i^\star,j}^{K-1}U_{i^\star,j}^{K-1}$, using the same selection and output rule as CBO-MPPI. Each final evaluated center is shifted and jittered independently for the next control period, rather than replacing the population by the selected output. An all-infeasible local batch has $\widehat F_i^k=+\infty$ and uses its projected center as its local update or output.

\begin{algorithm}[t]
\algocfconts{alg:paper_multistart}{\caption{MS-MPPI: one control period}}{
\begin{algorithmic}[1]
\REQUIRE State $x$; initialized or warm-started centers $M_i^0\in\R^{mH}$, $i=1,\ldots,N$; $K,L$, $\tau>0$, $\Sigma\succ0$; $\sigma_{\rm jit},B_{\rm jit}$
\ENSURE Applied input and next-period centers
\FOR{$k=0,\ldots,K-1$}
\FOR{all centers $i=1,\ldots,N$ in parallel}
\STATE Draw $\epsilon_{i,j}^k\sim\mathcal N(0,\Sigma)$; set $U_{i,j}^k=M_i^k+\epsilon_{i,j}^k$ and evaluate $S_{i,j}^k=\bar J_x(U_{i,j}^k)$ from $x$
\STATE Compute $\widehat F_i^k$ and the local weighted mean from~\eqref{eq:paper_multistart_weights}, with the zero-weight safeguard
\IF{$k<K-1$}
\STATE Update $M_i^{k+1}$ by~\eqref{eq:paper_multistart_update}
\ENDIF
\ENDFOR
\ENDFOR
\STATE Select $i^\star\in\argmin_i\widehat F_i^{K-1}$; reuse its final local mean as $U_{\rm out}$ and apply $[U_{\rm out}]_0$
\FOR{all $i=1,\ldots,N$}
\STATE Draw $A_i$ with independent $\mathcal N(0,1)$ entries; set $\zeta_i=\sigma_{\rm jit}B_{\rm jit}A_i$
\STATE Set $M_i^{0,\rm next}=\operatorname{Shift}(M_i^{K-1})+\zeta_i$
\ENDFOR
\end{algorithmic}}
\end{algorithm}
\subsubsection{CBO-MPC}
\label{subsec:cbo_mpc}
CBO-MPC treats each particle $U_i^k\in\R^D$ as a candidate control sequence and evaluates $S_i^k=\bar J_x(U_i^k)$ directly \citep{BorghiMonSep0100:00:00UTC2025Model,Sun2026Consensusbased}. For $\alpha>0$, its cost-based weights and consensus are
\begin{equation}
\nu_i^k=\frac{e^{-\alpha S_i^k}}{\sum_{r=1}^N e^{-\alpha S_r^k}},\qquad
U_\alpha^k=\sum_{i=1}^N\nu_i^k U_i^k.
\label{eq:cbo_mpc_consensus}
\end{equation}
After each of the first $K-1$ evaluations, the particles move according to
\begin{equation}
U_i^{k+1}=U_i^k-\lambda\Delta s(U_i^k-U_\alpha^k)
+\sigma_{\rm cbo}\|U_i^k-U_\alpha^k\|\xi_i^k\sqrt{\Delta s},
\quad \xi_i^k\sim\mathcal N(0,I_D).
\label{eq:cbo_mpc_discrete_update}
\end{equation}
Here $\lambda,\sigma_{\rm cbo},\Delta s>0$ have the same roles as in~\eqref{eq:discrete_cbo_mppi_update}. The final output is the minimum-cost particle: $i^\star\in\argmin_i S_i^{K-1}$ and $U_{\rm out}=U_{i^\star}^{K-1}$. There is no local MPPI batch or weighted-mean output. In the numerical implementation, an all-infeasible population uses the mean of its projected particles for consensus and output.

\begin{algorithm}[htbp]
\algocfconts{alg:cbo_mpc}{\caption{CBO-MPC: one control period}}{
\begin{algorithmic}[1]
\REQUIRE State $x$; initialized or warm-started particles $U_i^0\in\R^{mH}$, $i=1,\ldots,N$; $K$, $\alpha,\lambda,\sigma_{\rm cbo},\Delta s>0$; $\sigma_{\rm jit},B_{\rm jit}$
\ENSURE Applied input and next-period particles
\FOR{$k=0,\ldots,K-1$}
\STATE Roll out every $U_i^k$ from $x$ and evaluate $S_i^k=\bar J_x(U_i^k)$ in parallel
\STATE Compute $U_\alpha^k$ by~\eqref{eq:cbo_mpc_consensus}, using the all-infeasible safeguard if needed
\IF{$k<K-1$}
\STATE Draw independent $\xi_i^k\sim\mathcal N(0,I_D)$ and update particles by~\eqref{eq:cbo_mpc_discrete_update}
\ENDIF
\ENDFOR
\STATE Select $i^\star\in\argmin_i S_i^{K-1}$ and set $U_{\rm out}=U_{i^\star}^{K-1}$, using the all-infeasible safeguard if needed; apply $[U_{\rm out}]_0$
\FOR{all $i=1,\ldots,N$}
\STATE Draw $A_i$ with independent $\mathcal N(0,1)$ entries; set $\zeta_i=\sigma_{\rm jit}B_{\rm jit}A_i$
\STATE Set $U_i^{0,\rm next}=\operatorname{Shift}(U_i^{K-1})+\zeta_i$
\ENDFOR
\end{algorithmic}}
\end{algorithm}
\FloatBarrier

\section{Complete Numerical Results}
\label{sec:supplementary_numerics}
We give the common numerical protocol, task specifications, complete results,
and parameter-sensitivity analysis.
\subsection{Common Numerical Protocol}
\label{full:sec:numerical_example}
\label{sec:paper_experimental_details}
\label{sec:paper_method_details}

\paragraph{Common settings.}
Within each task, all four controllers use the same dynamics, trajectory cost, constraints, control period, prediction horizon $H$, closed-loop evaluation window $T$, input set, and total rollout budget $L_{\mathrm{tot}}$ per control step. These quantities vary across tasks and are specified in Appendices~\ref{subsec:paper_navigation_details}--\ref{subsec:paper_arm_details}. Preliminary timing experiments set each rollout budget so that CBO-MPPI's median run-level control-cycle 95th percentile was below the task's control period.

\paragraph{MPPI-based controller settings.}
MPPI, MS-MPPI, and CBO-MPPI use the same task-dependent temperature and local perturbation covariance. Preliminary free-energy landscapes determine the low- and high-temperature regimes, and preliminary small-scale experiments determine the reference MPPI and CBO-MPPI parameters. These task-dependent regimes give $\tau=3000$ for Many Obstacles and $\tau=300$ for Navigation. MPPI uses one nominal sequence and $L=L_{\mathrm{tot}}$ rollouts. MS-MPPI and CBO-MPPI share $(N,L,K)$, the initial-center law, local perturbations, center-selection and MPPI weighted-output rule, input projection, and warm start; they differ only in the center update, using~\eqref{eq:paper_multistart_update} and~\eqref{eq:discrete_cbo_mppi_update}, respectively.

CBO-MPPI uses $\lambda=1$, no gradient assistance, and the reference noise-to-drift ratio
\[
\kappa:=\frac{D\sigma_{\rm cbo}^2}{2\lambda}=0.5,
\qquad \sigma_{\rm cbo}=\sqrt{\frac{2\lambda\kappa}{D}}.
\]
The reference $\alpha$, $\Delta s$, covariance, initialization, warm-start jitter, and allocation vary across tasks and are given in the task sections.

\paragraph{CBO-MPC settings.}
\label{subsec:cbo_mpc_selection}
CBO-MPC inherits $\alpha$, $\lambda$, $\kappa$, $\sigma_{\rm cbo}$, $\Delta s$, the initialization law, warm start, input projection, and $K$ from CBO-MPPI. It sets $N_{\rm CBO\text{-}MPC}=N_{\rm CBO\text{-}MPPI}L_{\rm CBO\text{-}MPPI}$ to match the common rollout budget. CBO-MPC forms its cost-weighted consensus by~\eqref{eq:cbo_mpc_consensus}, whereas CBO-MPPI forms its free-energy-weighted consensus by~\eqref{eq:finite_sample_consensus}. The sensitivity study varies only $\alpha$ or $\kappa$ for CBO-MPPI; the horizon study varies $H$ and recomputes $\sigma_{\rm cbo}$ to preserve $\kappa$. All other quantities retain the corresponding task reference values.

\paragraph{Temporal correlation.}
This paragraph specifies the autoregressive correlation construction used for the sampling covariance $\Sigma$ in~\eqref{eq:mppi_free_energy} and its numerical counterparts: the local Gaussian perturbations in Algorithms~\ref{alg:cbo_mppi_free_energy}--\ref{alg:paper_multistart}, the initial center or particle populations in Algorithms~\ref{alg:cbo_mppi_free_energy}, \ref{alg:paper_multistart}, and~\ref{alg:cbo_mpc}, and the post-hoc evaluation perturbations in~\eqref{eq:paper_closed_loop_metrics}. It does not apply to the isotropic CBO diffusion vectors $\xi_i^k$ in~\eqref{eq:discrete_cbo_mppi_update} and~\eqref{eq:cbo_mpc_discrete_update}, or to the smooth warm-start jitter defined below. For controller sampling and initialization, write $\rho$ for the temporal correlation and $\sigma$ for the marginal standard deviation, and view an $H$-step perturbation as an array $\epsilon\in\mathbb R^{H\times m}$. Each channel $c$ is a stationary first-order autoregressive sequence,
\begin{equation}
  \epsilon_{0,c}=\sigma\eta_{0,c},
  \qquad
  \epsilon_{s,c}=\rho\,\epsilon_{s-1,c}+\sigma\sqrt{1-\rho^2}\,\eta_{s,c},
  \qquad s=1,\ldots,H-1,
  \label{eq:paper_ar1_perturbation}
\end{equation}
where the $\eta_{s,c}$ are independent standard Gaussian variables. Every time step therefore has marginal standard deviation $\sigma$ and $\operatorname{Cov}(\epsilon_{s,c},\epsilon_{t,c})=\sigma^2\rho^{|s-t|}$, so the sampling covariance $\Sigma$ of Section~\ref{sec:problem_setting} is $\sigma^2R_\rho\otimes I_m$ with $(R_\rho)_{st}=\rho^{|s-t|}$. We use $\rho=0.9$ for Navigation, Many Obstacles, and Pendulum. The 10-Hz Planar Arm discretization uses $\rho=0.81=0.9^2$, preserving the correlation over the same physical time separation. For the post-hoc evaluation in~\eqref{eq:paper_closed_loop_metrics}, we apply~\eqref{eq:paper_ar1_perturbation} over the full $T$-step evaluation window. The task sections give the corresponding $\sigma$ for each use.

\paragraph{Warm starts.}
Algorithms~\ref{alg:cbo_mppi_free_energy}, \ref{alg:paper_multistart}, and~\ref{alg:cbo_mpc} warm-start each center or particle by shifting it and adding the smooth jitter $\zeta_i$ defined below. We use $R_{\rm jit}=\min(4,H)$ cosine modes; each center or particle is viewed as an $H\times m$ array, and $A_i\in\R^{R_{\rm jit}\times m}$ contains independent standard Gaussian entries:
\begin{equation}
  \begin{gathered}
    \zeta_i=\sigma_{\rm jit}B_{\rm jit}A_i,
    \qquad
    \widetilde B_{sq}=\cos(\pi t_sq),
    \qquad t_s=\frac{s+\tfrac12}{H},\\
    s=0,\ldots,H-1,
    \qquad q=0,\ldots,R_{\rm jit}-1,
    \qquad
    B_{\rm jit}[s,:]=\widetilde B[s,:]\big/\lVert\widetilde B[s,:]\rVert_2 .
  \end{gathered}
  \label{eq:paper_jitter_basis}
\end{equation}
Because each row of $B_{\rm jit}$ is a unit vector, every time step receives jitter of marginal standard deviation $\sigma_{\rm jit}$. The basis restricts the jitter to the lowest-frequency modes, so it varies smoothly along the horizon. For the shift operations in Algorithms~\ref{alg:cbo_mppi_free_energy}, \ref{alg:paper_multistart}, and~\ref{alg:cbo_mpc}, the implementation appends zero.

\paragraph{Closed-loop evaluation.}
To evaluate $F_{\rm CL}$ in~\eqref{eq:paper_closed_loop_metrics}, we hold the realized input sequence fixed and roll out each entire perturbed sequence from $x_0$, without rerunning feedback control. Its cost and admissible set use the initial state $x_0$ and the full $T$-step evaluation window. Evaluation perturbations have the controller's marginal standard deviation and temporal correlation, extended from $H$ to $T$ steps. The same antithetic perturbations are reused across methods and trials within each task, with $L_{\rm e}=32768$ except for Many Obstacles, where $L_{\rm e}=8192$. The input-increment penalties specified in each task are included in $J_{\rm CL}$; state-dependent stage terms use $x_{k+1}$, which gives the convention of Section~\ref{sec:problem_setting} after composition with $f$.

Success criteria are given in each task section. In every task, Success is a first-hitting event within the $T$-step evaluation window: the trajectory need not remain in the target set afterward or belong to it at the terminal step $T$. Task-specific feasibility conditions are checked through the first hitting step. Arrival time $t_{\rm suc}$, in seconds, is summarized over successful trials only; minimum clearance $c_{\min}$ is measured over the full evaluation window. The task tables use MS-MPPI.

\subsection{Navigation}
\label{subsec:navigation}
\label{subsec:paper_navigation_details}
\label{full:subsec:two_obstacles}

\paragraph{Problem Setup.}
We use the point-mass instance of~\eqref{eq:discrete_time_dynamics} and the geometry specified in Section~\ref{subsec:two_obstacles}, with initial state $x_{\mathrm{init}}:=x_0$, initial position $p_{\mathrm{init}}:=p_0$, and target position $p_{\mathrm g}$.
The finite trajectory cost has the form
\begin{align}
  J_{x_0}(U)
  ={} &
  \sum_{k=0}^{H-1}
  \bigl(
  q_{\mathrm{goal}}\lVert p_{k+1}-p_{\mathrm g}\rVert^2
  +q_{\mathrm{vel}}\lVert v_{k+1}\rVert^2
  \bigr)
  +\sum_{k=0}^{H-1}
  \bigl(
  q_u\lVert u_k\rVert^2
  +q_{\Delta u}\lVert \Delta u_k\rVert^2
  \bigr)
  \notag                                                     \\
      & +q_{\mathrm{goal},f}\lVert p_H-p_{\mathrm g}\rVert^2
  +q_{\mathrm{vel},f}\lVert v_H\rVert^2,
  \label{eq:trajectory_cost}
\end{align}
where $\Delta u_0:=u_0-u_{\mathrm{prev}}$ and $\Delta u_k:=u_k-u_{k-1}$ for $k\geq1$, with $u_{\mathrm{prev}}$ denoting the control applied at the preceding control step and set to zero at the beginning of an episode. The cost is minimized over $U=(u_0,\ldots,u_{H-1})$ subject to a componentwise input bound $|u_{k,j}|\leq u_{\max}$.
The cost \eqref{eq:trajectory_cost} uses weights
\[
  q_{\mathrm{goal}}=0.1,\quad q_{\mathrm{vel}}=0.05,\quad q_u=0.01,
  \quad q_{\Delta u}=0.02,\quad q_{\mathrm{goal},f}=30,\quad q_{\mathrm{vel},f}=10,
\]
with input bound $u_{\max}=4$ and prediction horizon $H=60$. Closed-loop evaluation uses $T=80$ steps, corresponding to $8.0$~s.
For a realized trajectory, we define the minimum obstacle clearance over the evaluation window by
\[
  c_{\min}=\min_{0\leq k\leq T,o}
  (\lVert p_k-c_o\rVert-r_o),
\]
where $o$ indexes the obstacles and $c_o$ and $r_o$ denote the center and radius of obstacle $o$, respectively. Thus $c_{\min}=0$ is feasible contact and $c_{\min}<0$ is collision.

Under $L_{\mathrm{tot}}=131{,}072$, MPPI uses $L=131{,}072$, MS-MPPI and CBO-MPPI use $(N,L,K)=(256,64,8)$, and CBO-MPC uses $(N,K)=(16384,8)$. These allocations satisfy the common budget in Section~\ref{sec:numerical_example}.

For CBO-MPPI, we use $\alpha=0.02$, $\Delta s=0.2$ and $\sigma_{\mathrm{cbo}}=0.0913$, with initial-center standard deviation $2.0$.
At each control period, CBO-MPPI and CBO-MPC warm-start their shifted centers and particles, respectively, by adding independent smooth jitter with standard deviation $0.35$.
The local MPPI perturbation for CBO-MPPI has standard deviation $0.35$; within each CBO evaluation round, a shared perturbation is used across centers.

\paragraph{Open-Loop Free-Energy Landscape.}
\label{subsec:case1_landscape}

The routes in Figure~\ref{fig:paper_landscapes} (left) are generated from nominal sequences $\{M_h\}_{h\in[-3,3]}$ as follows.

Let $\phi_k=k/(H{-}1)$ denote normalized time and define $s(\phi)=\phi^2(1{-}\phi)^2$ and $d(\phi)=\phi^2(1{-}\phi)^2(1{-}2\phi)$.
We write the goal displacement as $\Delta p=p_{\mathrm g}{-}p_{\mathrm{init}}$ and the vertical unit vector as $e_y=(0,1)^\top$.
With normalizing constants $A_s$ and $A_h$, we set the velocity profile to
\begin{equation}
  v_k(h)=\frac{s(\phi_k)}{A_s}\,\Delta p
  +\frac{h}{A_h}\,d(\phi_k)\,e_y,
  \qquad k=0,\ldots,H-1,\qquad v_H(h)=0.
  \label{eq:route_height_velocity}
\end{equation}
The bell function $s$ is positive between the endpoints and vanishes at both endpoints. The antisymmetric function $d$ is positive on $(0,\tfrac12)$ and negative on $(\tfrac12,1)$.
We choose $A_s$ so that the position integral of the first term equals $\Delta p$, and $A_h$ so that the peak vertical deviation contributed by the second term equals~$h$.
The nominal control sequence $M_h=(u_0(h),\ldots,u_{H-1}(h))$ is obtained by finite-difference acceleration $u_k(h)=\Pi_{\mathcal C_u}\!\bigl((v_{k+1}(h){-}v_k(h))/\Delta t\bigr)$.

We use the following evaluation procedure for all landscapes below. Let $r$ denote the landscape horizon. We first project each parameterized nominal sequence onto $\mathcal C_u^r$. We then add perturbations to that projected nominal and project the resulting sequences again before rollout.

The nominal hard-constrained cost is $\bar J(h):=\bar J_{x_{\mathrm{init}}}(M_h)$.
Using $L_{\mathrm{MC}}$ control perturbations, the free-energy landscape estimator around $M_h$ is
\begin{equation}
  \widehat F(h;\tau)
  = -\tau\log\Bigl(\frac{1}{L_{\mathrm{MC}}}\sum_{\ell=1}^{L_{\mathrm{MC}}}
  e^{-\bar J_{x_{\mathrm{init}}}(\Pi_{\mathcal C_u^H}(M_h+\epsilon_\ell))/\tau}\Bigr),
  \label{eq:case1_landscape_F_estimator}
\end{equation}
where $L_{\mathrm{MC}}$ is the number of rollouts used only for this offline Monte Carlo landscape estimate, distinct from the controller batch size $L$ and the post-hoc evaluation count $L_{\mathrm e}$. The perturbations $\epsilon_\ell$ follow~\eqref{eq:paper_ar1_perturbation} with marginal standard deviation $0.35$ and temporal correlation $\rho=0.9$. Thus $\widehat F(h;\tau)$ is a Monte Carlo estimate of the projected free energy $F_x^\Pi(M_h)$.
Each curve in Figure~\ref{fig:paper_landscapes} (right) uses $L_{\mathrm{MC}}=65{,}536$ antithetic Monte Carlo rollouts.

In a finer sweep with $L_{\mathrm{MC}}=262{,}144$ antithetic rollouts, the transition temperature at the reference seed is approximately $21$, and the gap $\widehat F(0)-\min_h\widehat F(h)$ grows from about $6.0$ to $55.9$ and $215$ at $\tau=30$, $100$, and $300$.

\paragraph{Closed-Loop Evaluation Details.}
\label{subsec:case1_closed_loop_same_budget}

A navigation trial succeeds at the first step $k\leq T$ such that $\lVert p_k-p_{\mathrm g}\rVert\leq0.6$ and the entire realized trajectory through step $k$ is collision-free. Route classes use the crossing coordinate $y_{\mathrm{cross}}$: Center if $|y_{\mathrm{cross}}|\leq0.15$ and Detour otherwise.
Table~\ref{tab:case1_closed_loop_same_budget} adds Success, route counts, arrival times, clearance, and implementation timing to the cost comparison in Table~\ref{tab:paper_main_results}. The closed-loop trajectories of MPPI, CBO-MPC, and CBO-MPPI are shown in Figure~\ref{fig:paper_navigation}.

\begin{widepapertable}
  \centering
  \caption{Closed-loop performance for Navigation at $T=80$. Metrics are medians over 100 trials; arrival times use successful trials only. A dash in $t_{\rm suc}$ denotes no successful trial.}
  \label{tab:case1_closed_loop_same_budget}
  \small
  \setlength{\tabcolsep}{3pt}
  \begin{tabular}{lrrrrrrrrr}
    \hline
    Method   & $\tau$ & Success & Center & Detour & $t_{\mathrm{suc}}$ & $c_{\min}$ & $J_{\mathrm{CL}}$ & $F_{\mathrm{CL}}$ & $t_{\rm ctrl}$ [ms] \\
    \hline
    MPPI & 15 & 100/100 & 100 & 0 & 4.10 & 0.356 & 135.17 & 222.52 & 108.6 \\
    MS-MPPI & 15 & 99/100 & 70 & 30 & 3.10 & 0.253 & 163.05 & 224.64 & 48.3 \\
    CBO-MPC & 15 & 100/100 & 100 & 0 & 3.90 & 0.319 & 125.26 & 212.04 & 63.0 \\
    CBO-MPPI & 15 & 100/100 & 97 & 3 & 4.00 & 0.305 & 130.36 & 216.18 & 51.0 \\
    MPPI & 300 & 100/100 & 100 & 0 & 4.30 & 0.353 & 151.27 & 1015.66 & 108.6 \\
    MS-MPPI & 300 & 71/100 & 37 & 63 & 4.50 & 0.245 & 371.62 & 978.06 & 48.3 \\
    CBO-MPC & 300 & 100/100 & 100 & 0 & 3.90 & 0.319 & 125.26 & 955.57 & 63.0 \\
    CBO-MPPI & 300 & 100/100 & 6 & 94 & 5.30 & 0.717 & 151.98 & 820.67 & 51.0 \\
    \hline
  \end{tabular}
\end{widepapertable}

MS-MPPI succeeds in $99/100$ trials at $\tau=15$ and $71/100$ at $\tau=300$; the other three methods succeed in all trials. MS-MPPI's successful low-temperature trials arrive earlier than those of MPPI and CBO-MPPI. Arrival time is summarized over successful trials. At high temperature, CBO-MPPI's detours arrive later than MPPI's center routes.

\FloatBarrier
\subsection{Many Obstacles}
\label{subsec:many_obstacles}

\paragraph{Problem Setup.}
Many Obstacles uses the same point-mass instance of~\eqref{eq:discrete_time_dynamics} and hard-obstacle convention as Navigation, with fourteen circular obstacles over $[-10,10]^2$, initial state $(-8.5,-9,0,0)$, and goal $(8.5,9)$. The weights are $(q_{\mathrm{goal}},q_{\mathrm{vel}},q_u,\allowbreak q_{\Delta u},q_{\mathrm{goal},f},q_{\mathrm{vel},f})=\allowbreak(0.015,0.02,0.005,0.03,80,20)$, $u_{\max}=2.5$, and $H=50$. Closed-loop evaluation uses $T=120$ steps, corresponding to $12.0$~s. Under $L_{\mathrm{tot}}=131{,}072$, MPPI uses $L=131{,}072$, MS-MPPI and CBO-MPPI use $(N,L,K)=(256,64,8)$, and CBO-MPC uses $(N,K)=(16384,8)$. CBO-MPPI uses $\alpha=0.02$, $\sigma_{\mathrm{cbo}}=0.1$, and $\Delta s=0.2$, with initial-center standard deviation $2.5$ and warm-start jitter standard deviation $0.50$; the local perturbation has standard deviation $0.45$.

\paragraph{Closed-Loop Performance.}
A trial succeeds at the first step $k\leq T=120$ such that $\lVert p_k-p_{\mathrm g}\rVert\leq0.6$ and the entire realized trajectory through step $k$ is collision-free. To evaluate $F_{\rm CL}$, we use perturbation standard deviation $0.45$.

\begin{widepaperfigure}
  \centering
  \includegraphics[width=0.82\linewidth]{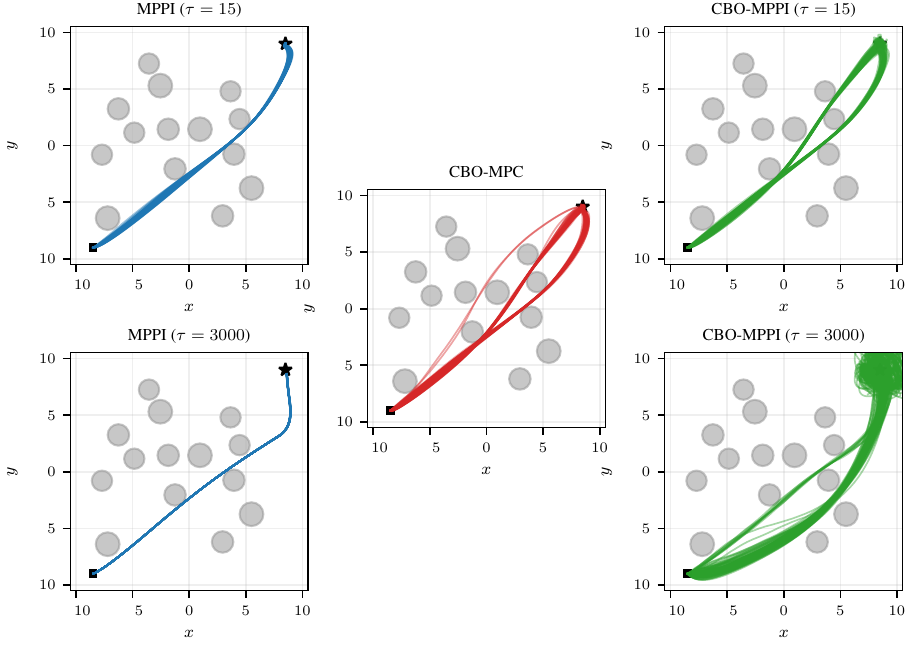}
  \caption{Five-panel comparison for Many Obstacles over the common 12-s window. The outer panels show MPPI and CBO-MPPI at $\tau=15$ (upper) and $\tau=3000$ (lower); the center shows CBO-MPC. Each panel overlays all 100 trials.}
  \label{fig:many_obstacles_closed_loop}
\end{widepaperfigure}

\begin{widepapertable}
  \centering
  \caption{Closed-loop performance for Many Obstacles at $T=120$. Metrics are medians over 100 trials; arrival times use successful trials only. A dash in $t_{\rm suc}$ denotes no successful trial.}
  \label{tab:many_obstacles_closed_loop_same_budget}
  \small
  \setlength{\tabcolsep}{4pt}
  \begin{tabular}{lrrrrrrr}
    \hline
    Method & $\tau$ & Success & $t_{\mathrm{suc}}$ & $c_{\min}$ & $J_{\mathrm{CL}}$ & $F_{\mathrm{CL}}$ & $t_{\rm ctrl}$ [ms] \\
    \hline
    MPPI & 15 & 100/100 & 5.80 & 0.068 & 243.10 & 424.24 & 120.3 \\
    MS-MPPI & 15 & 81/100 & 6.20 & 0.090 & 596.45 & 446.54 & 71.1 \\
    CBO-MPC & 15 & 100/100 & 6.10 & 0.080 & 239.39 & 417.21 & 83.4 \\
    CBO-MPPI & 15 & 100/100 & 6.20 & 0.084 & 250.69 & 430.88 & 72.8 \\
    MPPI & 3000 & 100/100 & 11.10 & 0.101 & 265.06 & 13708.01 & 120.3 \\
    MS-MPPI & 3000 & 40/100 & 7.60 & 0.453 & 1471.25 & 8120.80 & 71.1 \\
    CBO-MPC & 3000 & 100/100 & 6.10 & 0.080 & 239.39 & 14811.65 & 83.4 \\
    CBO-MPPI & 3000 & 90/100 & 6.00 & 0.458 & 456.85 & 7579.74 & 72.8 \\
    \hline
  \end{tabular}
\end{widepapertable}

At $\tau=15$, MPPI, CBO-MPC, and CBO-MPPI succeed in all 100 trials; MS-MPPI succeeds in 81. CBO-MPC has the lowest free energy, and MS-MPPI has the highest realized cost and free energy. At $\tau=3000$, CBO-MPPI succeeds in 90 trials, MS-MPPI in 40, and MPPI and CBO-MPC in all 100. Both MS-MPPI and CBO-MPPI lower free energy relative to MPPI and have higher realized cost. CBO-MPPI improves both costs and the success rate relative to MS-MPPI. Its larger-clearance trajectories reach the goal earlier among successful trials than MPPI and CBO-MPC.

\paperfloatbarrier

\subsection{Pendulum}
\label{subsec:pendulum}
\label{subsec:paper_pendulum_details}

We next conduct torque-limited pendulum swing-up experiments to test whether the $J$--$F_x^\Pi$ gap extends beyond hard-obstacle navigation.

\paragraph{Problem Setup.}

We impose the bound $|u|\leq u_{\max}$ on the actuator torque.
The state is $x=(\theta,\omega)$ with $\theta=0$ the upright target and $\theta=\pi$ the downward initial position.
The continuous-time dynamics are
\[
  \dot\theta = \omega, \qquad
  \dot\omega = \frac{g}{l}\sin\theta - \frac{c}{m_{\mathrm p}l^2}\omega + \frac{u}{m_{\mathrm p}l^2},
\]
discretized by the fourth-order Runge--Kutta method (RK4) with $\Delta t=0.05$, $g=9.81$, pendulum length $l=1$, mass $m_{\mathrm p}=1$, damping $c=0.05$, and input bound $u_{\max}=3$; the discrete model saturates the angular velocity at $|\omega|\leq 8$.
At the specified parameters, $u_{\max}/(m_{\mathrm p}gl)\approx 0.31$. The available torque is insufficient to reach the upright in a single half-swing. Swing-up therefore requires energy pumping over multiple oscillations.

Define the wrapped angle error by
\[
  e_\theta(\theta)
  :=
  \operatorname{atan2}(\sin\theta,\cos\theta).
\]
The regulation cost uses
\[
  \ell(x_t,u_t)
  =e_\theta(\theta_t)^2+0.1\,\omega_t^2,
  \qquad
  g(x_H)
  =e_\theta(\theta_H)^2+0.1\,\omega_H^2,
\]
with no control-penalty term.

A torque disturbance $\nu_t\sim\mathcal N(0,0.5^2)$ is added to each applied input.
The equal-budget comparison uses $L_{\mathrm{tot}}=16{,}384$ with planning horizon $H=80$. Closed-loop evaluation uses $T=170$ steps, corresponding to $8.5$~s.
MPPI uses $L=16{,}384$, MS-MPPI and CBO-MPPI use $(N,L,K)=(256,16,4)$, and CBO-MPC uses $(N,K)=(4096,4)$, so each method evaluates $16{,}384$ rollouts per step.
The MPPI-based controllers use perturbation standard deviation $0.25$. CBO-MPPI uses $\alpha=1$, $\sigma_{\mathrm{cbo}}=0.1118$, $\Delta s=0.1$, initial-center standard deviation $1.0$ with temporal correlation $0.9$, and warm-start jitter with standard deviation~$0.3$.
The closed-loop visualization and quantitative comparison use $\tau\in\{1,10\}$ for the MPPI-based controllers, together with the single temperature-independent CBO-MPC configuration.

\paragraph{Open-Loop Free-Energy Landscape.}
\label{subsec:pendulum_landscape}

We compare the swing-up strategies favored by $J$ and $F_x^\Pi$ in the obstacle-free Pendulum task. We use the landscape-evaluation procedure of Section~\ref{subsec:case1_landscape} with the following family of torque profiles.
Two parameters $(a,b)$ index the sinusoidal open-loop torque profiles
\[
  u_t(a,b)=a\sin\Bigl(\frac{2\pi\cdot 3\,t}{H'}\Bigr)
  +b\sin\Bigl(\frac{2\pi\cdot 4\,t}{H'}\Bigr),
  \qquad t=0,\ldots,H'-1.
\]
For this landscape calculation only, we use the longer horizon $H'=160$ (i.e., $8$\,s) with the same stage cost $\ell$ as the closed loop.
For each of the $81\times81$ grid points in $(a,b)$-space, $\widehat F(a,b;\tau)$ is estimated as in \eqref{eq:case1_landscape_F_estimator} from $L_{\mathrm{MC}}=4096$ perturbations following~\eqref{eq:paper_ar1_perturbation} with marginal standard deviation $0.25$ and temporal correlation $\rho=0.9$. Since this example has no state constraints, $\bar J_x$ equals $J_x$ on every input-projected sequence.

Free energy favors a swing-up strategy with a wider perturbation basin in this task.
Figure~\ref{fig:pendulum_landscape} shows the open-loop cost landscape $J(a,b)$ and the free-energy landscape $\widehat F(a,b;10)$ over this sinusoidal template family, with markers identifying the minimizer of each objective. To compare the trajectories' energy levels, let $E_{\mathrm{up}}=2m_{\mathrm p}gl$ denote the mechanical energy at the upright equilibrium, measured relative to the downward equilibrium. The minimizer of $J$ is a cost-efficient ``tightrope'' trajectory: it pumps the mechanical energy to $107\%$ of $E_{\mathrm{up}}$, crosses $\theta=0$ with nonzero angular velocity, and lies in a narrow valley of the cost landscape. The minimizer of $\widehat F$ uses a ``soft-landing'' strategy, pumping to almost exactly $100\%$ of $E_{\mathrm{up}}$ with near-zero angular velocity at the upright. This strategy has slightly higher nominal $J$ and lies in a wider low-cost neighborhood.

\begin{widepaperfigure}
  \centering
  \includegraphics[width=\linewidth]{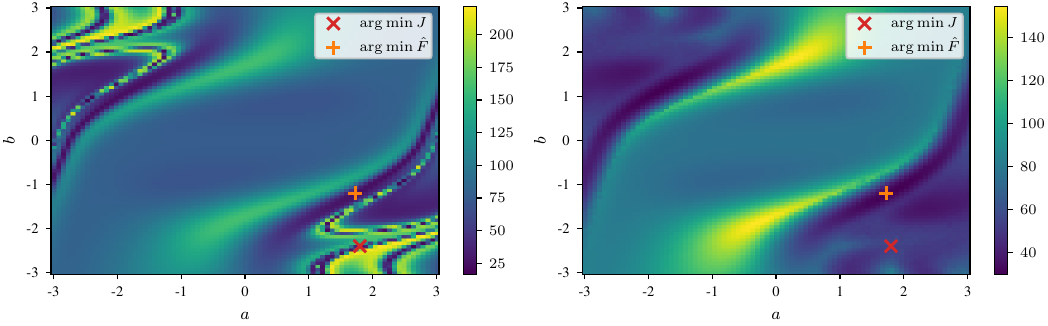}
  \caption{Open-loop landscape for Pendulum over sinusoidal templates ($H'=160$, $\tau=10$).
    Left: nominal cost $J(a,b)$. Right: free energy $\widehat F(a,b;10)$.
    Red~$\times$: $\arg\min J$ (energy overshoot $7\%$).
    Orange~$+$: $\arg\min \widehat F$ (overshoot $\le 0.2\%$).}
  \label{fig:pendulum_landscape}
\end{widepaperfigure}

\paragraph{Closed-Loop Performance.}
\label{subsec:pendulum_closed_loop}

We next compare swing-up attainment and closed-loop performance to assess the effect of the search mechanism under the same rollout budget.
A trial succeeds at the first state satisfying
\[
  |e_\theta(\theta_k)|<0.25,
  \qquad |\omega_k|<1.0,
  \qquad k\leq T.
\]
The realized cost $J_{\rm CL}$ is accumulated over the full evaluation window using the disturbed trajectory. To evaluate $F_{\rm CL}$, we use perturbation standard deviation $0.25$.
Figure~\ref{fig:pendulum_closed_loop} shows the phase portraits, and Table~\ref{tab:pendulum_closed_loop} summarizes the quantitative results.

\begin{widepaperfigure}
  \centering
  \includegraphics[width=\linewidth]{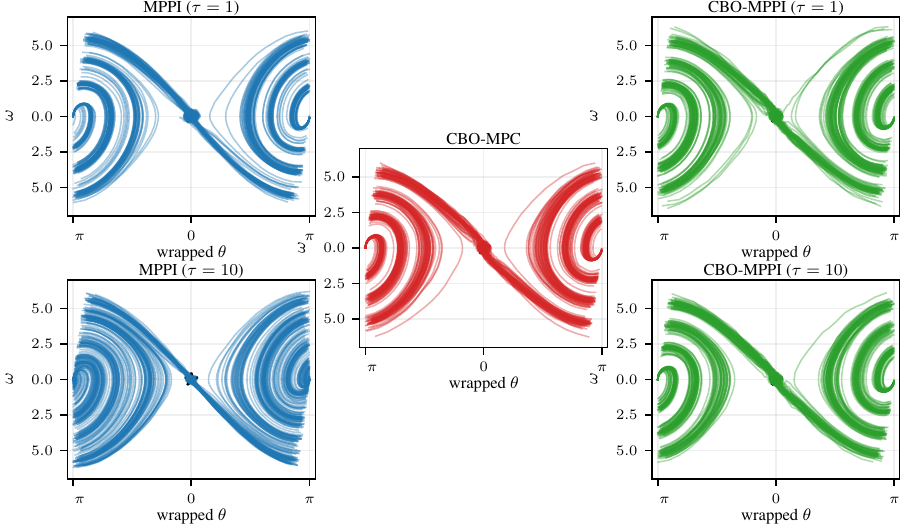}
  \caption{Five-panel phase-portrait comparison for Pendulum over the common $8.5$-s window. The outer panels show MPPI and CBO-MPPI at $\tau=1$ (upper) and $\tau=10$ (lower); the center shows the temperature-independent CBO-MPC baseline. Each panel overlays all 100 trials.}
  \label{fig:pendulum_closed_loop}
\end{widepaperfigure}

\begin{widepapertable}
  \centering
  \caption{Closed-loop performance for Pendulum at $T=170$. Metrics are medians over 100 trials; arrival times use successful trials only. A dash in $t_{\rm suc}$ denotes no successful trial.}
  \label{tab:pendulum_closed_loop}
  \small
  \setlength{\tabcolsep}{5pt}
  \begin{tabular}{lrrrrrr}
    \hline
    Method   & $\tau$ & Success & $t_{\mathrm{suc}}$ & $J_{\mathrm{CL}}$ & $F_{\mathrm{CL}}$ & $t_{\rm ctrl}$ [ms] \\
    \hline
    MPPI & 1 & 100/100 & 4.40 & 554.05 & 639.20 & 39.8 \\
    MS-MPPI & 1 & 96/100 & 5.95 & 688.35 & 679.30 & 34.9 \\
    CBO-MPC & 1 & 100/100 & 4.45 & 548.53 & 638.33 & 36.4 \\
    CBO-MPPI & 1 & 100/100 & 4.55 & 558.30 & 606.81 & 37.6 \\
    MPPI & 10 & 100/100 & 5.95 & 756.96 & 817.84 & 39.8 \\
    MS-MPPI & 10 & 98/100 & 6.30 & 743.83 & 738.78 & 34.9 \\
    CBO-MPC & 10 & 100/100 & 4.45 & 548.53 & 722.23 & 36.4 \\
    CBO-MPPI & 10 & 100/100 & 4.65 & 562.10 & 709.60 & 37.6 \\
    \hline
  \end{tabular}
\end{widepapertable}

MPPI, CBO-MPC, and CBO-MPPI succeed in all trials at both temperatures. MS-MPPI succeeds in $96/100$ trials at $\tau=1$ and $98/100$ at $\tau=10$. At $\tau=1$, MS-MPPI reaches the target later and has higher realized cost and free energy than MPPI and CBO-MPPI. CBO-MPPI has lower median free energy than MPPI together with slightly later arrival and higher realized cost. Appendix~\ref{full:subsec:parameter_sensitivity} shows the opposite low-temperature paired-median ordering throughout the parameter sweeps. At $\tau=10$, MS-MPPI improves both costs relative to MPPI, and CBO-MPPI improves them further and arrives earlier than MPPI. CBO-MPC has the lowest realized cost; CBO-MPPI has the lowest free energy.

\subsection{Planar Arm}
\label{subsec:planar_arm}
\label{subsec:paper_arm_details}

Finally, we conduct experiments with a redundant planar manipulator to investigate whether the free-energy objective favors larger-clearance postures by exploiting kinematic redundancy under hard link-obstacle constraints.

\paragraph{Problem Setup.}

The arm has three revolute joints with state $x=(q,v)\in\mathbb R^6$, where $q\in\mathbb R^3$ is the joint angle and $v\in\mathbb R^3$ is the joint velocity.
We use a joint-space damped double-integrator model, where the link lengths and trigonometric dependence enter through forward kinematics and clearance evaluation.

With time step $\Delta t=0.1$, damping $d=1$, and joint-space acceleration input $u\in\mathbb R^3$, the discrete dynamics are
\[
  q_{k+1}=q_k+\Delta t\,v_k,
  \qquad
  v_{k+1}=(1-d\Delta t)v_k+\Delta t\,u_k,
\]
with componentwise input bound $|u_{k,j}|\leq 5$.
For link lengths $\ell_1=\ell_2=\ell_3=1$, the joint positions are
\begin{equation}
  p_i(q)
  =
  \sum_{j=1}^{i}
  \ell_j
  \begin{bmatrix}
    \cos\left(\sum_{m=1}^{j}q_m\right) \\
    \sin\left(\sum_{m=1}^{j}q_m\right)
  \end{bmatrix},
  \qquad
  i=1,2,3,
  \label{eq:arm_forward_kinematics}
\end{equation}
with base point $p_0=(0,0)$ and end-effector position $p_{\mathrm{ee}}(q)=p_3(q)$.
The objective includes end-effector tracking, velocity damping, input effort, and input smoothness.
For an input sequence $U=(u_0,\ldots,u_{H-1})$, define the finite tracking and control cost
\begin{align}
  J_{x_0}(U)
  ={} &
  \sum_{k=0}^{H-1}
  \Bigl(
  q_{\mathrm{ee}}\lVert p_{\mathrm{ee}}(q_{k+1})-p_{\mathrm g}\rVert^2
  +q_v\lVert v_{k+1}\rVert^2
  \Bigr)
  +\sum_{k=0}^{H-1}
  \Bigl(
  q_u\lVert u_k\rVert^2
  +q_{\Delta u}\lVert \Delta u_k\rVert^2
  \Bigr)
  \notag                                                                  \\
      & +Q_{\mathrm{ee}}\lVert p_{\mathrm{ee}}(q_H)-p_{\mathrm g}\rVert^2
  +Q_v\lVert v_H\rVert^2,
  \label{eq:arm_finite_cost}
\end{align}
where $\Delta u_0:=u_0-u_{\mathrm{prev}}$ and $\Delta u_k:=u_k-u_{k-1}$ for $k\geq1$, with $u_{\mathrm{prev}}$ denoting the control applied at the preceding control step and set to zero at the beginning of an episode.
The numerical weights are $q_{\mathrm{ee}}=10$, $q_v=0.2$, $q_u=2\times10^{-3}$, $q_{\Delta u}=10^{-2}$, $Q_{\mathrm{ee}}=200$, and $Q_v=1$.
The prediction horizon is $H=20$, and closed-loop evaluation uses $T=35$ steps, corresponding to $3.5$~s.
The initial configuration is $q_0=(-0.8,1.0,-0.6)$, and the target end-effector position is $(1.8,1.2)$.
Two circular obstacles are placed at $(1.0,0.6)$ and $(1.5,0.9)$ with radii $0.25$ and $0.20$. We measure the link-obstacle clearance at configuration $q$ by
\begin{equation}
  c_{\mathrm{arm}}(q)
  =
  \min_{i,o}
  \left(
  d_{\mathrm{seg}}\bigl([p_{i-1}(q),p_i(q)],c_o\bigr)
  -r_o-r_{\ell}
  \right),
  \label{eq:arm_clearance}
\end{equation}
where $i$ and $o$ index the links and obstacles, respectively, $c_o$ and $r_o$ are the obstacle center and radius, $r_{\ell}=0.05$ is the link radius, and $d_{\mathrm{seg}}$ is the point-to-segment distance. For a realized trajectory, we define the minimum link-obstacle clearance over the evaluation window by
\[
  c_{\min}:=\min_{0\leq k\leq T}c_{\mathrm{arm}}(q_k).
\]
Zero clearance is feasible contact; negative clearance indicates link-obstacle penetration. After input projection, we enforce obstacle avoidance by assigning cost $+\infty$ to rollouts with $\min_{0\le k\le H}c_{\mathrm{arm}}(q_k)<0$. The same $\bar J_x$ is used in the landscape calculation, the closed-loop controllers, and all closed-loop metrics involving hard link-obstacle feasibility.
For the closed-loop comparison, the MPPI-based controllers are evaluated at $\tau\in\{3,300\}$ and use perturbation standard deviation $0.3$.
CBO-MPPI uses $\alpha=1$, $\sigma_{\mathrm{cbo}}=0.1291$, $\Delta s=0.1$, initial-center standard deviation $1.0$ with temporal correlation $0.81$, warm-start jitter standard deviation $0.2$, and $(N,L,K)=(256,64,8)$. MPPI uses $L=131{,}072$ rollouts per MPC step. CBO-MPC uses $(N,K)=(16384,8)$, $\alpha=1$, and $\sigma_{\rm cbo}=0.1291$; its $16384\times8$ evaluations give the same $L_{\mathrm{tot}}=131{,}072$ as MPPI and CBO-MPPI.

\paragraph{Open-Loop Free-Energy Landscape.}
\label{subsec:arm_landscape}

We evaluate $\bar J_x$ and $F_x^\Pi$ along the arm's redundant degree of freedom using the landscape-evaluation procedure of Section~\ref{subsec:case1_landscape}.
The three-link arm is kinematically redundant for a two-dimensional end-effector target. For this landscape calculation only, we use the 10-Hz discretization $\Delta t=0.1$~s and landscape horizon $H'=60$.
We construct a one-parameter family of nominal control sequences $M_h$ from a fixed closed-loop input sequence $U^{\mathrm{base}}=(u_0^{\mathrm{base}},\ldots,u_{H'-1}^{\mathrm{base}})$.
Let $q_\star$ be the final configuration reached by this base sequence, and let $n(q_\star)$ be a unit vector in the null space of the end-effector Jacobian at $q_\star$.
For normalized time $\phi_k=k/(H'-1)$, define the smooth ramp $r(\phi)=\phi^2(3-2\phi)$ and set $M_h=(u_0(h),\ldots,u_{H'-1}(h))$, where
\begin{equation}
  u_k(h)=\Pi_{\mathcal C_u}
  \left(u_k^{\mathrm{base}}+2h\,r(\phi_k)n(q_\star)\right),
  \qquad k=0,\ldots,H'-1,\quad h\in[-0.3,0.3].
  \label{eq:arm_posture_parameterization}
\end{equation}
The scalar $h$ therefore moves the terminal posture in the local null-space direction: the end-effector remains near the target as the elbow and wrist posture change, altering link-obstacle clearance.
For each $h$, the nominal hard-constrained objective is $\bar J(h)=\bar J_{x_0}(M_h)$.
Using $L_{\mathrm{MC}}=1024$ Gaussian perturbations with standard deviation $0.03$ and temporal correlation $\rho=0.81$, we estimate $F_x^\Pi$ by $\widehat F(h;\tau)$ as in \eqref{eq:case1_landscape_F_estimator}.
Figure~\ref{fig:arm_posture_landscape} shows representative final postures induced by several values of $h$ and the resulting $\bar J_x$ and $\widehat F$ landscapes.

\begin{widepaperfigure}
  \centering
  \newlength{\armPostureWidth}
  \newlength{\armLandscapeWidth}
  \newlength{\armPanelHeight}
  \newlength{\armPostureHeight}
  \newlength{\armLandscapeTopSpace}
  \setlength{\armPostureWidth}{0.4588\linewidth}
  \setlength{\armLandscapeWidth}{0.5065\linewidth}
  \setlength{\armPanelHeight}{0.382\linewidth}
  \setlength{\armPostureHeight}{0.3175\linewidth}
  \setlength{\armLandscapeTopSpace}{0.0196\linewidth}
  \subfigure[{Representative final postures.\label{fig:arm_posture_sweep}}]{
    \begin{minipage}[t]{\armPostureWidth}
      \vbox to \armPanelHeight{\offinterlineskip
        \hbox to \linewidth{\hfil\includegraphics[width=\linewidth,height=\armPostureHeight]{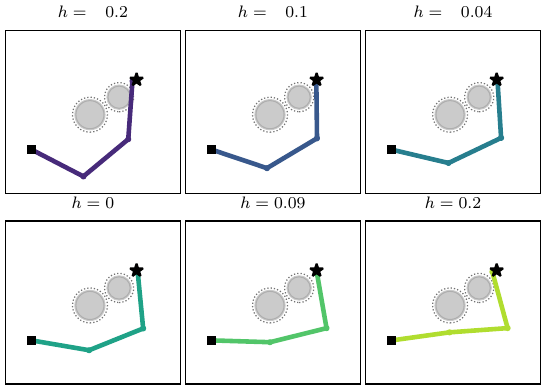}\hfil}
        \vfil}
    \end{minipage}}
  \hfill
  \subfigure[{Normalized $\bar J(h)$ and $\widehat F(h;\tau)$.\label{fig:arm_landscape}}]{
    \begin{minipage}[t]{\armLandscapeWidth}
      \vbox to \armPanelHeight{\offinterlineskip
        \kern\armLandscapeTopSpace
        \hbox to \linewidth{\hfil\includegraphics[width=\linewidth]{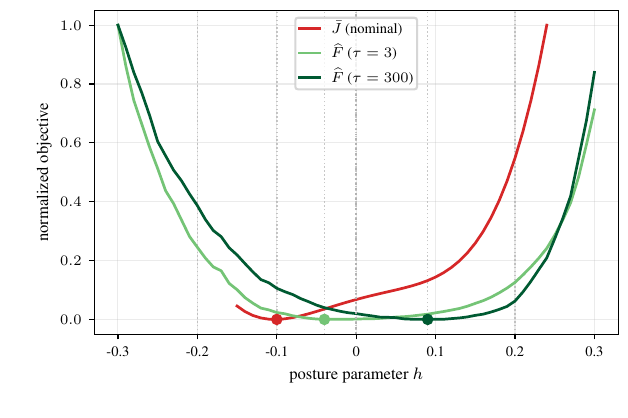}\hfil}
        \vfil}
    \end{minipage}}
  \caption{Posture parameterization and landscape for Planar Arm at 10~Hz ($H'=60$). Varying $h$ changes the redundant posture and keeps the end effector near the goal. Dotted circles in (a) have radius $r_o+r_{\ell}$ and mark the contact boundaries for the link centerlines, accounting for the link radius $r_{\ell}=0.05$. Curves in (b) are normalized to $[0,1]$.}
  \label{fig:arm_posture_landscape}
\end{widepaperfigure}

The preferred redundant posture moves toward larger-clearance neighborhoods under the free-energy objective, and the shift grows with temperature (Figure~\ref{fig:arm_posture_landscape}b). In the reference sweep, $\bar J_x$ is minimized at $h=-0.10$ with $c_{\mathrm{arm}}(q_{H'})=0.0276$. The free-energy minimizers move to $h=-0.04$ and $h=+0.09$ for $\tau=3$ and $\tau=300$, with $c_{\mathrm{arm}}(q_{H'})=0.0278$ and $0.0284$, respectively. The posture sweep uses perturbation standard deviation $0.03$; the closed-loop MPC comparison uses horizon $H=20$ and standard deviation $0.3$.

\paragraph{Closed-Loop Performance.}
\label{subsec:arm_closed_loop}

We next examine how objective choice and search mechanism affect target attainment and obstacle clearance in closed-loop operation.
A trial succeeds at the first step $k\leq T$ satisfying
\[
  \lVert p_{\mathrm{ee}}(q_k)-p_{\mathrm g}\rVert<0.20,
  \qquad \lVert v_k\rVert<1.5,
\]
provided that the entire realized trajectory through step $k$ is feasible, i.e., $c_{\mathrm{arm}}(q_j)\geq0$ for every $0\leq j\leq k$.
The realized cost $J_{\rm CL}$ uses \eqref{eq:arm_finite_cost} over the full evaluation window, and we use perturbation standard deviation $0.3$ to evaluate $F_{\rm CL}$.
Figure~\ref{fig:arm_closed_loop_final_postures} shows the final postures, and Table~\ref{tab:arm_closed_loop} summarizes the quantitative results.

\begin{widepaperfigure}
  \centering
  \includegraphics[width=\linewidth]{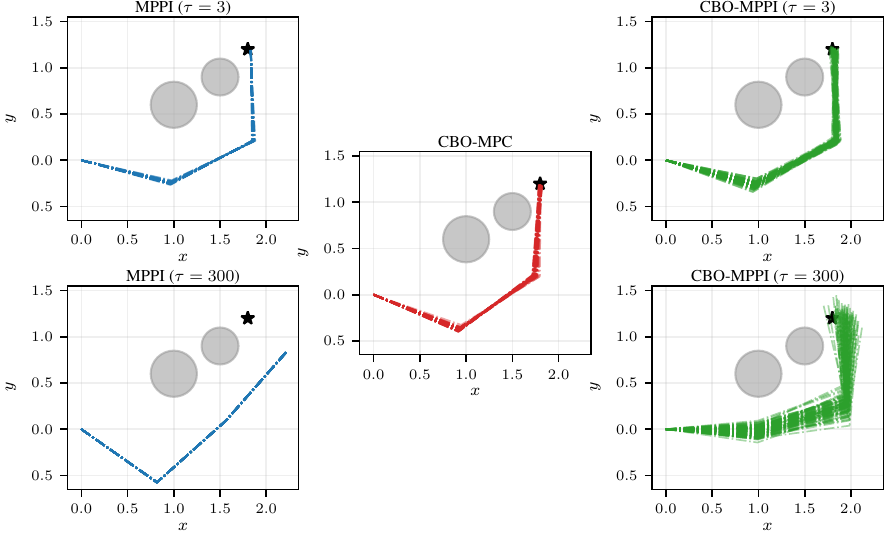}
  \caption{Five-panel comparison of final postures for Planar Arm at the end of the common $3.5$-s window. The outer panels show MPPI and CBO-MPPI at $\tau=3$ (upper) and $\tau=300$ (lower); the center shows the temperature-independent CBO-MPC baseline. Each panel overlays all 100 trials.}
  \label{fig:arm_closed_loop_final_postures}
\end{widepaperfigure}

\begin{widepapertable}
  \centering
  \caption{Closed-loop performance for Planar Arm at $T=35$. Metrics are medians over 100 trials; arrival times use successful trials only. A dash in $t_{\rm suc}$ denotes no successful trial.}
  \label{tab:arm_closed_loop}
  \small
  \setlength{\tabcolsep}{5pt}
  \begin{tabular}{lrrrrrrr}
    \hline
    Method & $\tau$ & Success & $t_{\mathrm{suc}}$ & $c_{\min}$ & $J_{\mathrm{CL}}$ & $F_{\mathrm{CL}}$ & $t_{\rm ctrl}$ [ms] \\
    \hline
    MPPI & 3 & 100/100 & 1.50 & 0.043 & 273.10 & 285.36 & 67.0 \\
    MS-MPPI & 3 & 100/100 & 1.00 & 0.024 & 184.77 & 205.71 & 90.0 \\
    CBO-MPC & 3 & 100/100 & 1.10 & 0.003 & 189.96 & 217.43 & 81.4 \\
    CBO-MPPI & 3 & 100/100 & 1.20 & 0.023 & 204.63 & 227.20 & 92.1 \\
    MPPI & 300 & 0/100 & -- & 0.330 & 694.06 & 1123.56 & 67.0 \\
    MS-MPPI & 300 & 71/100 & 2.60 & 0.164 & 332.81 & 789.72 & 90.0 \\
    CBO-MPC & 300 & 100/100 & 1.10 & 0.003 & 189.96 & 837.36 & 81.4 \\
    CBO-MPPI & 300 & 91/100 & 2.50 & 0.132 & 288.58 & 762.19 & 92.1 \\
    \hline
  \end{tabular}
\end{widepapertable}

At $\tau=3$, all four controllers succeed in all trials. MS-MPPI has the lowest realized cost and free energy and the earliest median arrival. CBO-MPC is second on both costs, while CBO-MPPI improves both over MPPI; independent local search remains best here. At $\tau=300$, MPPI never meets the target criterion, MS-MPPI succeeds in $71/100$ trials, CBO-MPPI succeeds in $91/100$, and CBO-MPC succeeds in all trials. CBO-MPC's $100/100$ count means that every trajectory enters the target set while feasible, with median first arrival at $1.10$~s; it does not require residence in the target set at the terminal time. Its direct optimization of trajectory cost aligns with the target-error and velocity conditions in the Success criterion, whereas high-temperature CBO-MPPI optimizes the perturbation-averaged objective and attains larger clearance with later arrival. CBO-MPPI has lower free energy and realized cost than MS-MPPI; MS-MPPI's larger median clearance accompanies a lower success rate. Relative to CBO-MPC, CBO-MPPI attains larger clearance ($0.132$ versus $0.003$), lower free energy, and higher realized cost.

\subsection{Effect of the Prediction Horizon}
\label{full:subsec:horizon_effect}

The reference settings fix one prediction horizon per task. Changing the horizon changes both the available lookahead and the dimension of the decision variable.

\paragraph{Protocol.}
For each task, we rerun the closed-loop experiment over a grid of horizons at that task's own $L_{\mathrm{tot}}$ and temperature, with twenty seeds per cell and the seed family of the task tables. The control period and the $T$-step evaluation window stay fixed as $H$ varies, so the panels compare controllers on one fixed control problem. We compare MPPI, MS-MPPI, and CBO-MPPI. MS-MPPI and CBO-MPPI share the allocation $(N,L,K)$, local perturbations, temperature, center selection, and warm start. MS-MPPI applies an MPPI update to each center; CBO-MPPI updates all centers through consensus drift and diffusion. CBO-MPPI keeps the noise-to-drift ratio $\kappa=D\sigma_{\rm cbo}^2/(2\lambda)$ of its reference setting, that is $\sigma_{\rm cbo}(H)=\sqrt{2\lambda\kappa/(mH)}$ with $m$ the input dimension, so the isotropic condition $2\lambda>D\sigma_{\rm cbo}^2$ holds at every horizon shown; MPPI keeps its reference perturbation scale. MPPI and CBO-MPPI optimize free energy, so we score each episode by both its realized cost $J_{\mathrm{CL}}$ and the free energy $F_{\mathrm{CL}}$ of the realized input sequence, as defined in~\eqref{eq:paper_closed_loop_metrics}. We compute $F_{\mathrm{CL}}$ post hoc at the higher of the two temperatures used for that task. The evaluation perturbation has fixed dimension $mT$, whereas a per-solve free energy has dimension $mH$.

\begin{widepaperfigure}
  \centering
  \includegraphics[width=\linewidth]{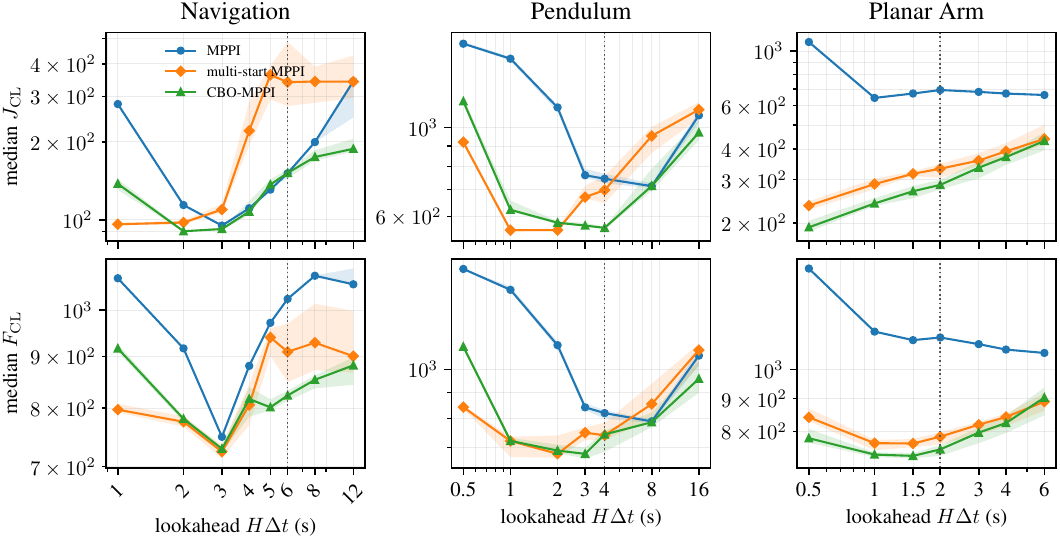}
  \caption{Realized cost (top) and free energy of the realized input sequence (bottom) against the prediction horizon, with the horizontal axis showing the lookahead $H\Delta t$. Each task uses its own $L_{\mathrm{tot}}$ ($131{,}072$ for Navigation and Planar Arm and $16{,}384$ for Pendulum) and reference temperature ($\tau=300$ for Navigation and Planar Arm and $\tau=10$ for Pendulum). All three controllers retain their other reference settings, CBO-MPPI keeping its reference $\kappa$; MS-MPPI differs from CBO-MPPI in the center update alone, so the vertical gap between the orange and green curves is the effect of the consensus term. Curves are medians over twenty seeds and bands are the interquartile range. The Planar Arm panel starts at the shortest tested horizon, $H=5$. The dotted line marks the horizon used elsewhere in the paper.}
  \label{fig:horizon_sweep_cost}
\end{widepaperfigure}

\paragraph{The preferred horizon depends on the task, controller, and criterion.}
The curves in Figure~\ref{fig:horizon_sweep_cost} are generally nonmonotone, and their minima span an order of magnitude across tasks. Free energy is lowest at $H=30$ on Navigation for all three controllers; on Pendulum, the minima occur at $H=160$, $40$, and $60$ for MPPI, MS-MPPI, and CBO-MPPI; and on Planar Arm, they occur at $H=60$, $15$, and $15$. The corresponding realized-cost minima occur at $H=(30,10,20)$ on Navigation, $H=(160,20,80)$ on Pendulum, and $H=(10,5,5)$ on Planar Arm, in the same controller order. Thus realized cost favors the shortest tested horizon for MS-MPPI and CBO-MPPI on Planar Arm, with medians $235.4$ and $192.1$. Longer horizons can also carry a substantial penalty: on Pendulum, moving from $H=80$ to $H=320$ raises median realized cost by $44\%$, $59\%$, and $73\%$ for MPPI, MS-MPPI, and CBO-MPPI, and leaves them with $6/20$, $1/20$, and $16/20$ successes.

\paragraph{The consensus term often sustains the long-horizon regime.}
MS-MPPI and CBO-MPPI draw the same number of rollouts from the same number of centers, allowing a direct comparison of their center updates. At short horizons, independent centers suffice and sometimes perform better: at $H=10$ on Navigation, MS-MPPI attains the lower realized cost ($96.2$ against $137.8$) and free energy ($797.1$ against $916.6$), and at $H=10$ on Pendulum, it succeeds in $15/20$ trials where CBO-MPPI succeeds in none, with both the lower cost ($919.9$ against $1163.8$) and the lower free energy ($841.3$ against $1109.2$). As the horizon grows, CBO-MPPI retains more successes on Navigation and Pendulum: MS-MPPI falls to $15/20$ at $H=60$ and $11/20$ at $H=120$ on Navigation, and to $7/20$ at $H=160$ and $1/20$ at $H=320$ on Pendulum, whereas CBO-MPPI records $19/20$, $17/20$, $20/20$, and $16/20$. On Planar Arm, CBO-MPPI leads $14/20$ to $5/20$ at $H=30$, while the two methods have similar success at $H=40$ and equal success at $H=60$. At the three reference horizons, CBO-MPPI has lower realized cost and higher success than MS-MPPI: $151.7$ against $340.1$ and $19/20$ against $15/20$ on Navigation; $561.6$ against $697.8$ and $20/20$ against $19/20$ on Pendulum; and $285.8$ against $332.4$ and $16/20$ against $15/20$ on Planar Arm. Its free energy is also lower on Navigation ($823.4$ against $909.5$) and Planar Arm ($750.1$ against $785.4$), while MS-MPPI is slightly lower on Pendulum ($740.0$ against $742.7$). At $H=60$ on Planar Arm, MS-MPPI similarly has lower free energy ($891.0$ against $904.1$), CBO-MPPI has slightly lower realized cost ($431.1$ against $438.9$), and both succeed in $7/20$ trials. MPPI trails MS-MPPI and CBO-MPPI over most of the grid: it meets the Planar Arm success criterion at no horizon, succeeds in $6/20$ trials at $H=320$ on Pendulum, and falls to $6/20$ at $H=120$ on Navigation.

\subsection{Parameter Sensitivity}
\label{full:subsec:parameter_sensitivity}

Using the reference controller settings in Appendices~\ref{sec:paper_experimental_details}--\ref{subsec:paper_arm_details}, we vary one parameter at a time over the grids
\[
  \alpha/\alpha_{\mathrm{ref}}\in\{0.25,0.5,1,2,4,8\},
  \qquad \kappa=D\sigma_{\mathrm{cbo}}^2/(2\lambda)\in\{0.1,0.3,0.5,0.7,0.9\}.
\]
Each parameter setting uses the same 20 independent trials paired seed-wise with MPPI. For trial $r$ within a fixed task--temperature case, let $F_{\mathrm{CBO}}^{(r)}(\theta)$ denote the CBO-MPPI $F_{\mathrm{CL}}$ at a swept parameter value $\theta$, and let $F_{\mathrm{MPPI}}^{(r)}$ denote the paired MPPI $F_{\mathrm{CL}}$. We measure sensitivity to the reference CBO-MPPI setting and the difference from MPPI using the following medians of paired differences over trials $r$:
\[
  \begin{aligned}
    \Delta F_{\mathrm{ref}}(\theta)
     & :=\operatorname*{median}_r\!
    \left[F_{\mathrm{CBO}}^{(r)}(\theta)-F_{\mathrm{CBO}}^{(r)}(\theta_{\mathrm{ref}})\right], \\
    \Delta F_{\mathrm{MPPI}}(\theta)
     & :=\operatorname*{median}_r\!
    \left[F_{\mathrm{CBO}}^{(r)}(\theta)-F_{\mathrm{MPPI}}^{(r)}\right].
  \end{aligned}
\]
Here $\theta_{\mathrm{ref}}$ is $\alpha/\alpha_{\mathrm{ref}}=1$ or $\kappa=0.5$. Negative values indicate lower free energy than the reference CBO-MPPI setting and paired MPPI, respectively. Because the median is taken after forming each seed-wise difference, these quantities may differ from differences of the corresponding direct medians. Tables~\ref{tab:sensitivity_alpha}--\ref{tab:sensitivity_kappa} report all settings and lower-$F$ counts.

Across the 72 non-reference settings, $\Delta F_{\mathrm{ref}}$ is negative at 34 and positive at 38, and every one of the 16 task--temperature--parameter sweeps contains at least one non-reference setting with negative $\Delta F_{\mathrm{ref}}$. Thus non-reference settings both increase and decrease the achieved free energy relative to the reference setting. After identifying the reference setting shared by the two parameter sweeps in each task--temperature case, the experiment contains 80 unique settings. Relative to MPPI, $\Delta F_{\mathrm{MPPI}}<0$ at 62 of these 80 settings: all 40 high-temperature settings and 22 of 40 low-temperature settings. Pendulum at $\tau=1$ is the exception, with $\Delta F_{\mathrm{MPPI}}>0$ throughout both sweeps.

In the tables, $F_{\rm MPPI}$ and $F_{\rm CBO\text{-}MPPI}$ denote median $F_{\rm CL}$. The lower-$F$ counts complement the paired medians by showing how often CBO-MPPI has smaller $F_{\rm CL}$ than paired MPPI.

\begin{widepapertable}
\centering
\caption{Complete $\alpha/\alpha_{\rm ref}$ sensitivity. Direct $F$ columns are medians; difference columns are medians of seed-wise paired differences and may differ from differences of the displayed direct medians. Lower-$F$ counts trials with $F_{\rm CBO\text{-}MPPI}<F_{\rm MPPI}$. Each parameter setting summarizes 20 paired trials.}
\label{tab:sensitivity_alpha}
\scriptsize
\setlength{\tabcolsep}{2.8pt}
\begin{tabular*}{\textwidth}{@{\extracolsep{\fill}}llrrrrrrr@{}}
\hline
Task & $\tau$ & $\alpha/\alpha_{\rm ref}$ & $F_{\rm MPPI}$ & $F_{\rm CBO\text{-}MPPI}$ & $\Delta F_{\rm ref}$ & $\Delta F_{\rm MPPI}$ & lower-$F$ & Success \\
\hline
Navigation & $15$ & $0.25$ & $227.2$ & $233.6$ & $10.7$ & $6.4$ & 4/20 & 20/20 \\
 &  & $0.5$ & $227.2$ & $227.8$ & $6.9$ & $0.7$ & 9/20 & 20/20 \\
 &  & $1$ & $227.2$ & $221.8$ & $0.0$ & $-5.4$ & 16/20 & 20/20 \\
 &  & $2$ & $227.2$ & $214.6$ & $-7.0$ & $-12.8$ & 19/20 & 20/20 \\
 &  & $4$ & $227.2$ & $207.6$ & $-14.5$ & $-19.0$ & 20/20 & 20/20 \\
 &  & $8$ & $227.2$ & $199.4$ & $-21.4$ & $-27.6$ & 20/20 & 20/20 \\
\cline{2-9}
 & $300$ & $0.25$ & $1045.6$ & $1021.8$ & $179.9$ & $-23.6$ & 17/20 & 20/20 \\
 &  & $0.5$ & $1045.6$ & $860.6$ & $18.9$ & $-184.8$ & 19/20 & 20/20 \\
 &  & $1$ & $1045.6$ & $833.0$ & $0.0$ & $-212.9$ & 20/20 & 20/20 \\
 &  & $2$ & $1045.6$ & $821.0$ & $-14.5$ & $-224.6$ & 20/20 & 20/20 \\
 &  & $4$ & $1045.6$ & $808.3$ & $-24.9$ & $-237.2$ & 20/20 & 20/20 \\
 &  & $8$ & $1045.6$ & $776.0$ & $-52.8$ & $-269.5$ & 20/20 & 20/20 \\
\hline
Pendulum & $1$ & $0.25$ & $642.4$ & $668.0$ & $-5.9$ & $40.2$ & 6/20 & 20/20 \\
 &  & $0.5$ & $642.4$ & $662.7$ & $3.0$ & $24.7$ & 9/20 & 20/20 \\
 &  & $1$ & $642.4$ & $659.4$ & $0.0$ & $6.1$ & 9/20 & 20/20 \\
 &  & $2$ & $642.4$ & $669.2$ & $0.5$ & $6.9$ & 9/20 & 20/20 \\
 &  & $4$ & $642.4$ & $667.9$ & $-5.3$ & $14.3$ & 7/20 & 20/20 \\
 &  & $8$ & $642.4$ & $673.3$ & $-4.0$ & $4.9$ & 9/20 & 20/20 \\
\cline{2-9}
 & $10$ & $0.25$ & $797.6$ & $688.9$ & $-12.5$ & $-96.2$ & 20/20 & 20/20 \\
 &  & $0.5$ & $797.6$ & $735.3$ & $3.1$ & $-74.7$ & 19/20 & 20/20 \\
 &  & $1$ & $797.6$ & $732.1$ & $0.0$ & $-83.8$ & 19/20 & 20/20 \\
 &  & $2$ & $797.6$ & $732.2$ & $-9.3$ & $-83.8$ & 19/20 & 20/20 \\
 &  & $4$ & $797.6$ & $732.7$ & $-2.6$ & $-74.5$ & 20/20 & 20/20 \\
 &  & $8$ & $797.6$ & $734.7$ & $-2.9$ & $-65.9$ & 19/20 & 20/20 \\
\hline
Planar Arm & $3$ & $0.25$ & $289.6$ & $243.7$ & $17.4$ & $-44.3$ & 20/20 & 20/20 \\
 &  & $0.5$ & $289.6$ & $239.7$ & $7.8$ & $-51.4$ & 20/20 & 20/20 \\
 &  & $1$ & $289.6$ & $225.7$ & $0.0$ & $-61.8$ & 20/20 & 20/20 \\
 &  & $2$ & $289.6$ & $228.4$ & $-1.9$ & $-60.0$ & 20/20 & 20/20 \\
 &  & $4$ & $289.6$ & $223.2$ & $-4.0$ & $-65.2$ & 20/20 & 20/20 \\
 &  & $8$ & $289.6$ & $224.3$ & $-4.3$ & $-64.7$ & 20/20 & 20/20 \\
\cline{2-9}
 & $300$ & $0.25$ & $1124.3$ & $756.2$ & $-5.7$ & $-368.6$ & 20/20 & 17/20 \\
 &  & $0.5$ & $1124.3$ & $760.8$ & $1.6$ & $-364.0$ & 20/20 & 18/20 \\
 &  & $1$ & $1124.3$ & $760.5$ & $0.0$ & $-364.4$ & 20/20 & 18/20 \\
 &  & $2$ & $1124.3$ & $770.4$ & $11.6$ & $-354.3$ & 20/20 & 17/20 \\
 &  & $4$ & $1124.3$ & $761.5$ & $6.3$ & $-362.7$ & 20/20 & 18/20 \\
 &  & $8$ & $1124.3$ & $766.5$ & $9.4$ & $-358.5$ & 20/20 & 15/20 \\
\hline
Many Obstacles & $15$ & $0.25$ & $424.6$ & $436.5$ & $-6.9$ & $9.2$ & 8/20 & 20/20 \\
 &  & $0.5$ & $424.6$ & $434.3$ & $1.8$ & $8.7$ & 9/20 & 20/20 \\
 &  & $1$ & $424.6$ & $430.1$ & $0.0$ & $5.2$ & 7/20 & 20/20 \\
 &  & $2$ & $424.6$ & $432.5$ & $2.2$ & $4.2$ & 9/20 & 20/20 \\
 &  & $4$ & $424.6$ & $417.1$ & $-17.5$ & $-14.7$ & 13/20 & 20/20 \\
 &  & $8$ & $424.6$ & $416.2$ & $-14.0$ & $-9.3$ & 14/20 & 20/20 \\
\cline{2-9}
 & $3000$ & $0.25$ & $13790.4$ & $8731.4$ & $1077.3$ & $-5068.8$ & 20/20 & 20/20 \\
 &  & $0.5$ & $13790.4$ & $8106.0$ & $363.8$ & $-5695.1$ & 20/20 & 20/20 \\
 &  & $1$ & $13790.4$ & $7632.3$ & $0.0$ & $-6160.6$ & 20/20 & 17/20 \\
 &  & $2$ & $13790.4$ & $7456.7$ & $-163.4$ & $-6328.8$ & 20/20 & 18/20 \\
 &  & $4$ & $13790.4$ & $7204.8$ & $-418.5$ & $-6574.9$ & 20/20 & 20/20 \\
 &  & $8$ & $13790.4$ & $7214.2$ & $-234.6$ & $-6576.4$ & 20/20 & 20/20 \\
\hline
\end{tabular*}
\end{widepapertable}

\clearpage

\begin{widepapertable}
\centering
\caption{Complete $\kappa$ sensitivity. Direct $F$ columns are medians; difference columns are medians of seed-wise paired differences and may differ from differences of the displayed direct medians. Lower-$F$ counts trials with $F_{\rm CBO\text{-}MPPI}<F_{\rm MPPI}$.}
\label{tab:sensitivity_kappa}
\scriptsize
\setlength{\tabcolsep}{2.8pt}
\begin{tabular*}{\textwidth}{@{\extracolsep{\fill}}llrrrrrrr@{}}
\hline
Task & $\tau$ & $\kappa$ & $F_{\rm MPPI}$ & $F_{\rm CBO\text{-}MPPI}$ & $\Delta F_{\rm ref}$ & $\Delta F_{\rm MPPI}$ & lower-$F$ & Success \\
\hline
Navigation & $15$ & $0.1$ & $227.2$ & $225.2$ & $3.4$ & $-2.0$ & 14/20 & 20/20 \\
 &  & $0.3$ & $227.2$ & $224.4$ & $2.4$ & $-2.8$ & 14/20 & 20/20 \\
 &  & $0.5$ & $227.2$ & $221.8$ & $0.0$ & $-5.4$ & 16/20 & 20/20 \\
 &  & $0.7$ & $227.2$ & $217.5$ & $-6.2$ & $-9.6$ & 19/20 & 20/20 \\
 &  & $0.9$ & $227.2$ & $216.6$ & $-4.1$ & $-10.7$ & 18/20 & 20/20 \\
\cline{2-9}
 & $300$ & $0.1$ & $1045.6$ & $839.3$ & $1.4$ & $-206.6$ & 20/20 & 20/20 \\
 &  & $0.3$ & $1045.6$ & $839.9$ & $3.3$ & $-205.9$ & 20/20 & 20/20 \\
 &  & $0.5$ & $1045.6$ & $833.0$ & $0.0$ & $-212.9$ & 20/20 & 20/20 \\
 &  & $0.7$ & $1045.6$ & $823.3$ & $-10.2$ & $-222.5$ & 20/20 & 20/20 \\
 &  & $0.9$ & $1045.6$ & $722.6$ & $-102.7$ & $-323.1$ & 20/20 & 5/20 \\
\hline
Pendulum & $1$ & $0.1$ & $642.4$ & $664.7$ & $2.3$ & $12.1$ & 7/20 & 20/20 \\
 &  & $0.3$ & $642.4$ & $627.8$ & $-2.6$ & $14.0$ & 6/20 & 20/20 \\
 &  & $0.5$ & $642.4$ & $659.4$ & $0.0$ & $6.1$ & 9/20 & 20/20 \\
 &  & $0.7$ & $642.4$ & $664.8$ & $0.2$ & $6.8$ & 8/20 & 20/20 \\
 &  & $0.9$ & $642.4$ & $679.9$ & $27.0$ & $15.5$ & 6/20 & 20/20 \\
\cline{2-9}
 & $10$ & $0.1$ & $797.6$ & $734.5$ & $5.1$ & $-86.1$ & 19/20 & 20/20 \\
 &  & $0.3$ & $797.6$ & $742.6$ & $4.0$ & $-65.4$ & 19/20 & 20/20 \\
 &  & $0.5$ & $797.6$ & $732.1$ & $0.0$ & $-83.8$ & 19/20 & 20/20 \\
 &  & $0.7$ & $797.6$ & $713.6$ & $-0.7$ & $-91.2$ & 20/20 & 20/20 \\
 &  & $0.9$ & $797.6$ & $737.1$ & $13.1$ & $-57.5$ & 18/20 & 20/20 \\
\hline
Planar Arm & $3$ & $0.1$ & $289.6$ & $274.5$ & $46.0$ & $-14.0$ & 19/20 & 20/20 \\
 &  & $0.3$ & $289.6$ & $255.8$ & $27.8$ & $-31.6$ & 20/20 & 20/20 \\
 &  & $0.5$ & $289.6$ & $225.7$ & $0.0$ & $-61.8$ & 20/20 & 20/20 \\
 &  & $0.7$ & $289.6$ & $204.6$ & $-24.0$ & $-84.0$ & 20/20 & 20/20 \\
 &  & $0.9$ & $289.6$ & $261.3$ & $29.1$ & $-26.7$ & 17/20 & 20/20 \\
\cline{2-9}
 & $300$ & $0.1$ & $1124.3$ & $808.0$ & $40.5$ & $-317.0$ & 20/20 & 16/20 \\
 &  & $0.3$ & $1124.3$ & $784.0$ & $18.4$ & $-339.7$ & 20/20 & 15/20 \\
 &  & $0.5$ & $1124.3$ & $760.5$ & $0.0$ & $-364.4$ & 20/20 & 18/20 \\
 &  & $0.7$ & $1124.3$ & $768.7$ & $14.1$ & $-355.3$ & 20/20 & 17/20 \\
 &  & $0.9$ & $1124.3$ & $754.9$ & $-25.1$ & $-369.4$ & 20/20 & 13/20 \\
\hline
Many Obstacles & $15$ & $0.1$ & $424.6$ & $409.7$ & $-27.4$ & $-17.8$ & 13/20 & 20/20 \\
 &  & $0.3$ & $424.6$ & $439.1$ & $20.9$ & $14.8$ & 8/20 & 20/20 \\
 &  & $0.5$ & $424.6$ & $430.1$ & $0.0$ & $5.2$ & 7/20 & 20/20 \\
 &  & $0.7$ & $424.6$ & $427.5$ & $-6.0$ & $-0.3$ & 10/20 & 20/20 \\
 &  & $0.9$ & $424.6$ & $515.6$ & $76.3$ & $97.0$ & 1/20 & 7/20 \\
\cline{2-9}
 & $3000$ & $0.1$ & $13790.4$ & $8042.7$ & $425.2$ & $-5758.5$ & 20/20 & 19/20 \\
 &  & $0.3$ & $13790.4$ & $8158.9$ & $544.9$ & $-5643.9$ & 20/20 & 19/20 \\
 &  & $0.5$ & $13790.4$ & $7632.3$ & $0.0$ & $-6160.6$ & 20/20 & 17/20 \\
 &  & $0.7$ & $13790.4$ & $7145.1$ & $-593.4$ & $-6643.0$ & 20/20 & 19/20 \\
 &  & $0.9$ & $13790.4$ & $11486.6$ & $4081.9$ & $-2289.2$ & 15/20 & 7/20 \\
\hline
\end{tabular*}
\end{widepapertable}

Some parameter settings combine lower free energy with lower success rates. For example, Many Obstacles at $\tau=15$, $\kappa=0.9$ succeeds in $7/20$ trials, and Navigation at $\tau=300$, $\kappa=0.9$ succeeds in $5/20$. The Many Obstacles $\kappa=0.9$ cells retain three numerical failures at $\tau=15$ and four at $\tau=3000$ in the 20-seed denominators. Both cases lie in the tested range $\kappa<1$, equivalently $2\lambda-D\sigma_{\mathrm{cbo}}^2>0$, used by the idealized mean-field theory.